\documentclass[12pt]{amsart}
\usepackage[all]{xy}
\usepackage{amsmath,amssymb,dsfont,mathtools,enumitem}
\usepackage[pdftitle={K-regularity of locally convex algebras},
  pdfauthor={Hvedri Inassaridze},
  pdfsubject={Mathematics}
]{hyperref}

\usepackage[lite]{amsrefs}
		
\usepackage{microtype}
\usepackage{wasysym}

\newtheorem{thm}{Theorem}[section]
\newtheorem{cor}[thm]{Corollary}

\newtheorem{prop}[thm]{Proposition}
\theoremstyle{definition}
\newtheorem{defn}[thm]{Definition}
\theoremstyle{remark}
\newtheorem{rem}[thm]{Remark}

\numberwithin{equation}{section}

\begin{document}

\title{(Co)homology of $\Gamma$-groups and $\Gamma$-homological algebra}
\author{Hvedri Inassaridze}
\address{A.~Razmadze Mathematical Institute of Tbilisi State University, 6, Tamarashvili Str., Tbilisi 0179, Georgia.}
\email{inassari@gmail.com}

\thanks{}

\begin{abstract}
This is a further investigation of our approach to group actions in homological algebra in the settings of homology of $\Gamma$-simplicial groups, particularly of $\Gamma$-equivariant homology and cohomology of $\Gamma$-groups. This approach could be called $\Gamma$-homological algebra. The abstract kernel of non-abelian extensions of groups, its relation with the obstruction to the existence of non-abelian extensions and with the second group cohomology are extended to the case of non-abelian $\Gamma$-extensions of $\Gamma$-groups. We compute the rational $\Gamma$-equivariant (co)homology groups of finite cyclic $\Gamma$-groups. The isomorphism of the group of n-fold $\Gamma$-equivariant extensions of a $\Gamma$-group G by a $G\rtimes \Gamma$-module A with the (n+1)th $\Gamma$-equivariant group cohomology of G with coefficients in A is proven.We define the $\Gamma$-equivariant Hochschild homology as the homology of the $\Gamma$- Hochschild complex when the action of the group $\Gamma$ on the Hochschild complex is induced by its action on the basic ring. Important properties of the $\Gamma$-equivariant Hochschild homology related to Kahler differentials, Morita equivalence and derived functors are established. Group (co)homology and $\Gamma$-equivariant group (co)homology of crossed $\Gamma$-modules are introduced and investigated by using relevant derived functors. Relations with extensions of crossed $\Gamma$-modules, in particular with  relative extensions of group epimorphisms in the sense of Loday and with $\Gamma$-equivariant extensions of crossed $\Gamma$-modules are established. Universal and $\Gamma$-equivariant universal central $\Gamma$-extensions of $\Gamma$-perfect crossed $\Gamma$-modules are constructed and Hopf formulas for the integral homology and $\Gamma$-equivariant integral homology of crossed $\Gamma$-modules are obtained. Finally, applications to algebraic K-theory, Galois theory of commutative rings and cohomological dimension of groups are given.

\end{abstract}

\maketitle

\section{introduction}

We continue the study of our approach to group actions in homological algebra which we call $\Gamma$-homological algebra that was started in [27] and continued in [28]. The origin of the equivariant study of group extensions theory in homological algebra goes back to Whitehead paper [48]. Group actions on algebraic and topological objects have many important applications in K-theory and homotopy theory([6,16,31,41]). Our goal is to continue the development of extension theory in the category of $\Gamma$-groups and of the relevant equivariant (co)homology theory that has been initiated in [27,28]. A different (co)homology theory of groups with operators was provided and investigated in [7-10], motivated by the graded categorical groups classification problem [8]. The introduction of $\Gamma$-equivariant chain complexes and their homology groups substantially contribute to the realization of our aim. Moreover this approach allows us to present a version of equivariant Hochschild homology of any unital k-algebra A induced by the action of the group $\Gamma$ on the k-algebra A. The $\Gamma$-equivariant Hochschild homology is closely related to $\Gamma$- equivariant homology of groups [27]. This equivariant version differs of the equivariant Hochschild homology given in [38].

We study extensions of $\Gamma$-groups that can be viewed as a part of group actions in homological algebra, particularly of group actions on simplicial groups. Two important classes of $\Gamma$-group extensions are considered. The first class is consisting of extensions having $\Gamma$-section map. The investigation of these extensions was initiated in [27] and called $\Gamma$-equivariant extensions of $\Gamma$-groups by $\Gamma$-equivariant G-modules. For the second class we deal with extensions of $\Gamma$-groups endowed with a crossed $\Gamma$-module structure and called $\Gamma$-extensions of crossed $\Gamma$-modules (having $\Gamma$-section map). Our approach to extensions of crossed modules substantially extends the class of relative extensions of group epimorphisms introduced and investigated by Loday [33]. It should be noted that homology and cohomology of crossed modules related to extensions of crossed modules were investigated by many authors [3,7,12-13,15-17,23].

The study of $\Gamma$-group extensions having $\Gamma$-section map is closely related to the extension problem of group actions satisfying some conditions, in our case to lifting group actions that split for the given extension of groups.

Applications to algebraic K-theory, Galois theory of commutative rings and cohomological dimension of $\Gamma$-groups are given.

The paper is divided into seven sections:

2. Preliminaries

3. Extensions of $\Gamma$-groups,

4. Some computations,

5. $\Gamma$-derived functors and $\Gamma$-equivariant Hochschild homology,

6. Extensions of crossed $\Gamma$-modules,

7. Homology and central $\Gamma$-extensions of crossed $\Gamma$-modules,

8. Applications to algebraic K-theory, Galois theory of commutative rings and cohomological dimension of groups.

Some notation that will be used throughout the paper:

$[\Gamma G]$ denotes the set of elements $^{\gamma}gg^{-1}$ and $\Gamma G$ is the normal subgroup of G generated by the set $[\Gamma G]$.
$G_\Gamma$ denotes the quotient group $G/\Gamma G$.

$[G,G]_{\Gamma}$ denotes the subgroup of the $\Gamma$-group G generated by the commutant subgroup [G,G] and the elements of the form $^{\gamma}gg^{-1}$, $g\in G$, $\gamma \in \Gamma$. It is called the $\Gamma$-commutant of the group G.

$[G,H]_{\Gamma}$ denotes the subgroup of G generated by the elements $x^{\gamma}yx^{-1}y^{-1}$, where $g\in G$, y belongs to the normal subgroup H of G, and $\gamma \in \Gamma$.

$G^{ab}_\Gamma$ denotes the abelianization of the group $G_\Gamma$.

\;


\section{Preliminaries}

In this section we recall some definitions and propositions given in [27] which will be used later. Moreover it is shown how equivariant versions of well known homological properties of groups are obtained by using our approach to group actions in homological algebra.

Let ${\mathbf{G}}^\Gamma$ be the category whose objects are groups on which a fixed group $\Gamma$ is acting, called $\Gamma$-groups, and morphisms are group homomorphisms compatible with the action of $\Gamma$.

Any exact sequence E of $\Gamma$-groups
\begin{equation}
E: 1\rightarrow A\rightarrow B\overset{\tau}\rightarrow G\rightarrow 1
\end{equation}
is called $\Gamma$-extension of the $\Gamma$-group G by the $\Gamma$-group A. The extension E is an extension with $\Gamma$-section map if there is a map $\beta: G\rightarrow X$ such that $\tau\beta = 1_G$ and $\beta$ is compatible with the action of $\Gamma$, $\beta (^{\gamma}g) = ^\gamma{\beta (g), g\in G, \sigma\in \Gamma}$. In addition if $\beta$ is a homomorphism then the extension E is called split extension.
\begin{defn}
(1) A $\Gamma$-equivariant G-module A is a G-module equipped with a $\Gamma$-module structure and the actions of G and $\Gamma$ are related to each other by the equality
$$
^{\sigma}(^{g}a)=^{\sigma}g(^{\sigma} a)
$$
for $g\in G$, $\sigma\in \Gamma$,$a\in A$.

The category of $\Gamma$-equivariant G-modules is equivalent to the category of $G\rtimes \Gamma$-modules, where $G\rtimes \Gamma$ is the semi-direct product of G and $\Gamma$ [9].

If E is an extension with $\Gamma$-section map and A is a $\Gamma$-equivariant G-module it is called  $\Gamma$-equivariant extension of G by A. In addition if X and G are $\Gamma$-equivariant G-modules it is called proper sequence of $\Gamma$-equivariant G-modules.

(2) A $\Gamma$-equivariant G-module F is called relatively free $G\rtimes\Gamma$-module if it is a free G-module with basis a $\Gamma$-set and relatively projective $\Gamma$-equivariant G-modules are retracts of relatively free $\Gamma$-equivariant G-modules.
\end{defn}

The class $\mathcal{P}$ of relatively projective $\Gamma$-equivariant G-modules is a projective class with respect to proper sequences of $\Gamma$-equivariant G-modules.

For the cohomological description of the set $E^{1}_{\Gamma}(G,A)$ of equivalence classes of $\Gamma$-equivariant extensions of G by A the $\Gamma$-equivariant homology and cohomology of $\Gamma$-groups have been introduced as relative $Tor_n^{\mathcal{P}}$ and $Ext^n_{\mathcal{P}}, n \geq 0$ in the category of $\Gamma$-equivariant G-modules [27], namely
\begin{defn}
The $\Gamma$-equivariant homology and cohomology of $\Gamma$-groups are defined as follows
$$
 H_n^{\Gamma}(G,A)= Tor_n^{\mathcal{P}}(\mathbb{Z},A), H^n_{\Gamma}(G,A)= Ext^n_{\mathcal{P}}(\mathbb{Z},A), n\geq 0,
$$
where the functors $\bigotimes$ and $Hom$ are taken over the ring $\mathbb{Z}(G\rtimes \Gamma)$ and the groups G and $\Gamma$ are trivially acting on the abelian group $\mathbb{Z}$ of integers.
\end{defn}

Let
\begin{equation}
\cdot\cdot\cdot\rightarrow B_n\rightarrow \cdot\cdot\cdot\rightarrow B_1\rightarrow B_0\rightarrow \mathbb{Z}\rightarrow 0
\end{equation}
be the bar resolution of $\mathbb{Z}$, where $B_0 = \mathbb{Z}(G)$ and $B_n,n>0$,is the free  $\mathbb{Z}(G)$-module generated by $[g_1,g_2,...,g_n]$, $g_i\in G$. The action of $\Gamma$ on G induces an action of $\Gamma$ on the sequence (2.2) defined by $^{\gamma}(mg)= m^{\gamma}g$ and
 $^{\gamma}(g[g_1,g_2,...,g_n]) = ^{\gamma}{g}[^{\gamma}{g_1},^{\gamma}{g_2},...,^{\gamma}{g_n}]$ for $B_0$ and $B_n$ respectively, $n\geq 1$. Then the sequence 2.2 is the $\Gamma$-equivariant bar resolution of $\mathbb{Z}$, the groups $B_n$ being relatively free $\Gamma$-equivariant G-modules and there are isomorphisms $H_n^{\Gamma}(G,A) \cong H_n(B_*{\bigotimes}_{G\ltimes \Gamma} A)$, $H^n_{\Gamma}(G,A) \cong H_n(Hom_{G\ltimes \Gamma}(B_*,A))$.

 For the $\Gamma$-equivariant cohomology of $\Gamma$-groups an alternative description by cocycles is provided. To this end the group $C^{n}_{\Gamma}(G,A)$ of $\Gamma$-maps, $f:G^{n}\rightarrow A$ for $n> 0$, called $\Gamma$-cochains is considered. By using the classical cobord operators $\delta^{n}:C^{n}_{\Gamma}(G,A)\rightarrow C^{n+1}_{\Gamma}(G,A)$, $n > 0$, we obtain a cochain complex
 $$
 0\rightarrow C^{0}_{\Gamma}(G,A)\rightarrow C^{1}_{\Gamma}(G,A)\rightarrow C^{2}_{\Gamma}(G,A)\rightarrow \cdot \cdot \cdot\rightarrow C^{n}_{\Gamma}(G,A)\rightarrow \cdot \cdot \cdot,
 $$
 where $C^{0}_{\Gamma}(G,A) = A^{\Gamma}$, $Ker\delta^{1} = Der_{\Gamma}(G,A)$ is the group of $\Gamma$-derivations and the homology groups of the complex $C^{*}_{\Gamma}(G,A)$ are isomorphic to the $\Gamma$-equivariant cohomology groups of the $\Gamma$-group G with coefficients in the $\Gamma$-equivariant G-module A.

 Two $\Gamma$-equivariant extensions E and E' of G by A are called equivalent if there is a morphism $E\rightarrow E'$ which is the identity on G and A. We denote by $E^{1}_{\Gamma}(G,A)$ the set of equivalence classes of $\Gamma$-equivariant extensions of G by A.
 \begin{thm}
 There is a bijection
 $E^{1}_{\Gamma}(G,A)\cong  H^2_{\Gamma}(G,A).$
  \end{thm}

  \begin{rem}
 By using the Baer sum operation the set $E^{1}_{\Gamma}(G,A)$ becomes an abelian group and the bijection of Theorem 2.3 is an isomorphism. This theorem will be extended to higher dimensions for $n \geq 2$ by introducing the notion of n-fold $\Gamma$-equivariant extension of G by A (see Theorem 3.7).
 \end{rem}
 \begin{defn}
 A $\Gamma$-group is called $\Gamma$-free if it is a free group with basis a $\Gamma$-set.
 \end{defn}
 Any free group F(G) generated by a $\Gamma$-group G becomes a $\Gamma$-free group by the following action of $\Gamma$: $^\gamma{\mid g \mid} =  \mid ^\gamma{g} \mid $, $g\in G, \gamma \in \Gamma$. The defining property of the $\Gamma$-free group F with basis E is that every $\Gamma$-map $E\overset{f}\rightarrow G$ to a $\Gamma$-group G is uniquely extended to a $\Gamma$-homorphism $F\overset{f'}\rightarrow G$.

 Let $\mathbb{F}$ be the projective class of $\Gamma$-free groups in the category ${\mathbf{G}}^\Gamma$ of $\Gamma$-groups.
\begin{thm}
 There are isomorphisms
$$
H_n^{\Gamma}(G,A) \cong L^{\mathbb{F}}_{n-1}(I(G)\otimes_{G\ltimes \Gamma}A), H^n_{\Gamma}(G,A)\cong R^{n-1}_{\mathbb{F}}Der_{\Gamma}(G,A)
$$
for $n\geq 2$, where I(G) is the kernel of the natural homomorphism $\mathbb{Z}(G)\rightarrow \mathbb{Z}$ of $\Gamma$-equivariant G-modules, Der(G,A) is the group of $\Gamma$-derivations, $L^{\mathbb{F}}_{n-1}$ and R$^{n-1}_{\mathbb{F}}$ denote respectively the left and right derived functors with respect to the projective class $\mathbb{F}$.
\end{thm}

We also recall some results on $\Gamma$-equivariant integral homology [27]. These homology groups are simply denoted $H^{\Gamma}_n(G)$ for $H^{\Gamma}_n(G,\mathbb{Z})$, the groups G and $\Gamma$ acting trivially on $\mathbb{Z}$.
\begin{thm}
(1)There is an isomorphism
$$
L^{\mathbb{F}}_{n}(G^{ab}_\Gamma) \cong H^{\Gamma}_{n+1}(G)
$$
for $n\geq 0$.

(2) There are exact sequences
$$
0\rightarrow \Gamma G/[G,G]\cap \Gamma G\rightarrow H_1(G)\rightarrow H^{\Gamma}_1(G)\rightarrow 0,
$$
$$
\cdot\cdot\cdot\rightarrow H^{\Gamma}_{n+1}(G)\rightarrow L^{\mathbb{F}}_{n-1}U(G)\rightarrow H_n(G)\rightarrow H^{\Gamma}_n(G)
$$
$$
\rightarrow L^{\mathbb{F}}_{n-2}U(G)\rightarrow\nleftarrow \cdot\cdot\cdot\rightarrow H^{\Gamma}_3(G)\rightarrow L^{\mathbb{F}}_{1}U(G)\rightarrow H_2(G)
$$
$$
\rightarrow H^{\Gamma}_2(G)\rightarrow L^{\mathbb{F}}_{0}U(G)\rightarrow H_1(G)\rightarrow H^{\Gamma}_1(G)\rightarrow 0,
$$
relating $\Gamma$-equivariant integral homology with the classical integral homology of groups, where U is a covariant functor assigning to any $\Gamma$-group G the abelian group $[G,G]_{\Gamma}/[G,G]$.
\end{thm}
\begin{thm}
Let
$$
1\rightarrow A\rightarrow B\overset{\tau}\rightarrow G\rightarrow 1
$$
be a short exact sequence of $\Gamma$-groups with $\Gamma$-section map and $\alpha: P\rightarrow B$ be a $\Gamma$-projective presentation of the $\Gamma$-group B. Then there is an exact sequence
$$
0\rightarrow V\rightarrow H^{\Gamma}_2(B)\rightarrow H^{\Gamma}_2(G)\rightarrow A/[B,A]_{\Gamma}\rightarrow H^{\Gamma}_1(B)\rightarrow H^{\Gamma}_1(G)\rightarrow 0,
$$
where V is the kernel of the $\Gamma$-homomorphism $[P,S]_{\Gamma}/[P,R]_{\Gamma}\rightarrow [B,A]_{\Gamma}$ induced by $\alpha$, R = Ker $\alpha$ and S = Ker $\tau\alpha$.
\end{thm}
\begin{thm}
If G is a $\Gamma$-group, then

(1)
$H^{\Gamma}_1(G,A) = G/[G,G]_{\Gamma}\otimes A$,
G and $\Gamma$ are trivially acting on A.

(2)
$H^{\Gamma}_2(G)$ is isomorphic to the group $(R\cap [P,P]_{\Gamma})/[P,R]_{\Gamma}$, where $R = Ker \alpha$ and $\alpha: P\rightarrow G$ is a $\Gamma$-projective presentation of G (Hopf formula for the $\Gamma$-equivariant homology of groups).
\end{thm}

The Brown - Ellis formula is also obtained extending Hopf formula to higher $\Gamma$-equivariant homology of groups (see [28]).

\begin{defn}
A $\Gamma$-subgroup L of a $\Gamma$-group G is called retract of G if there is a $\Gamma$-homomorphism $f:G\rightarrow L$ sucht that its restriction to L is the identity map.
\end{defn}

\begin{thm}

Let L be retract of a $\Gamma$-free group F. For any $\Gamma$-equivariant G-module A we have

1) exact sequences

$0\rightarrow H^{\Gamma}_{1}(L,A)\rightarrow I(L)\otimes_{\mathbb{Z}(L\rtimes \Gamma)} A\rightarrow A_{\Gamma}\rightarrow A_{L\rtimes \Gamma}\rightarrow 0$ and

$0\rightarrow\rightarrow A^{L\rtimes \Gamma}\rightarrow A^{\Gamma}\rightarrow Hom_{\mathbb{Z}(L\rtimes\Gamma)}(I(L),A)\rightarrow H^{1}_{\Gamma}(L,A)\rightarrow 0$,

2) $H_{n}^{\Gamma}(L,A) = 0$ and $H^{n}_{\Gamma}(L,A) = 0$ for $n> 1$.

\end{thm}

The consideration of retracts of $\Gamma$-free groups has motivated the following

$Question$: Let G be a $\Gamma$-subgroup of a $\Gamma$-free group. What are reasonable conditions, automatically satisfied in the case of trivial G, under which G is a $\Gamma$-free group? In other words, we are asking whether the well-known Nielsen--Schreier theorem on free groups can be nicely extended to $\Gamma$-free groups.

As G. Janelidze informed me, at least some partial answers, including the counter-example and Theorem 2.12 below, should be known (although we could not find a proper reference):

We take $\Gamma$ = ( 1, $\gamma$) a two-element cyclic group;	F to be the free group on a two-element set, say $\{x,y\}\}$ and $f : F \rightarrow \Gamma$ to be the group homomorphism defined by f(x) = $\gamma$ = f(y); define G to be Ker f. We make F a $\Gamma$-group, defining the $\Gamma$-action on it by $^{\gamma}x = y$ and $^{\gamma}y = x.$	The group $\Gamma$ is a $\Gamma$-group, assuming that $\Gamma$ acts	on itself trivially; this makes f a $\Gamma$-group homomorphism, and so G is a $\Gamma$-subgroup of F. After that we observe:
(a)	G is a free group on a three element set. Indeed, according to Nielsen—Schreier formula,   we have rank(G) = (F : G)(rank(F) — 1) + 1 = 2(2 — 1) + 1 = 3.
(b)	For	$\gamma \in \Gamma$  and $t \in  F$ , we have	$^{\gamma}t= t \Rightarrow	t = 1$. Indeed, if t is a non-empty word,  the t and $^{\gamma}t$	must begin with different letters.
(c)	Let B  be any  basis of  G;  since it  is  a  basis,  1  doesn't belongs to  B.  If  B  is a  $\Gamma$-subset  of G, then, presenting it as the disjoint union of $\Gamma$-orbits, we obtain a contradiction between (a) and (b). Indeed,  (a)  implies  that  at  least  one orbit must have exactly one element, which contradicts to (b) since 1 doesn't belongs to B. That means G is not $\Gamma$-free group.

Since  Nielsen - Schreier theorem on free groups doesn't hold for $\Gamma$-free groups, there is a meaning to study the structure of $\Gamma$-subgroups of $\Gamma$-free groups, in particular to establish conditions for  $\Gamma$-subgroups of $\Gamma$-free groups to be $\Gamma$-free groups.

For instance, let F be a $\Gamma$-group and G a $\Gamma$-subgroup of F. The action of $\Gamma$ on F induces an action of $\Gamma$ on the set $G/F = \{Gt|t\in F\}$ of right cosets such that the canonical map $p:F\rightarrow G/F$ is a surjection of $\Gamma$-sets.

\begin{thm}(G.Janelidze)

Let F be a $\Gamma$-free group and G a $\Gamma$-subgroup of F. Suppose there exists a $\Gamma$-map $s:G/F\rightarrow F$ satisfying the following conditions:
(1) we have ps = 1, (2) there exists a $\Gamma$-subset X of F such that F is a free group on X and $y_{1}...y_{n}y_{n+1}\in s(G/F)\Rightarrow y_{1}...y_{n}\in s(G/F)$ whenever $y_{1}...y_{n}y_{n+1}$ is a canonical length n+1 presentation of an element of F with $y_{i}\in X\bigcup \{x^{-1}|x\in X\}$ for each i = 1,...,n+1. Then G is a $\Gamma$-free group.

\end{thm}

As an application of this theorem consider the following simple example: We take $\Gamma = \{1,\gamma,\gamma^{2}\}$ a three-element cyclic group, F to be the $\Gamma$-free group on the three-element set \{x,y,z\}and defining the $\Gamma$-action on it by $^{\gamma}x = y, ^{\gamma}y = z, ^{\gamma}z = x$. We also take $\mathds{Z}/2\mathds{Z} \times \mathds{Z}/2\mathds{Z}$ to be a $\Gamma$-group assuming that $\Gamma$ acts on it by $^{\gamma}(c,0)= (0,c)$ and $^{\gamma}(0,c) = (c,c)$ where c denotes the generator of $\mathds{Z}/2\mathds{Z}$, and $f:F\rightarrow  \mathds{Z}/2\mathds{Z} \times \mathds{Z}/2\mathds{Z}$ to be a $\Gamma$-group homomorphism defined by f(x) = (c,0), f(y) = (0,c) and f(z) = (c,c). Then $G = Ker(f)$ is a $\Gamma$-free group.

Some results related to this problem are obtained in [49,50].

Finally we provide an assertion establishing relation of the $\Gamma$-equivariant cohomology of groups with the well known equivariant cohomology of topological spaces.

Let G be a $\Gamma$-group and X a topological space on which the groups G and $\Gamma$ are acting such that G is acting properly and
$^{\gamma}(^{g}x) = ^{\gamma}g(^{\gamma}x)$. $\gamma\in \Gamma, g\in G, x\in X.$
\begin{thm}
If X is either acyclic and $\Gamma$ acts trivially on X or X is $\Gamma$-contractible, then there is an isomorphism

$ H^n_{\Gamma}(G,A) \cong H^n_{\Gamma}(X/G,A)$ for $n\geq 0$, where G and $\Gamma$ are trivially acting on the abelian group A and $H^*_{\Gamma}(X/G,A)$ is the equivariant cohomology of the space X/G.

\end{thm}
\;

\section{Extensions of $\Gamma$-groups}

We introduce an internal property of $\Gamma$-group extensions possessing a $\Gamma$-section map that will be used through out the paper.

\begin{defn}

It will be said that the sequence 2.1 of $\Gamma$-groups possesses the $\Gamma$-property if the restriction of $\tau$ on the subset $[\Gamma B]$ of the group B is injective.

\end{defn}
\begin{thm}

The sequence 2.1 possesses the $\Gamma$-property iff it has a $\Gamma$-section map and $\Gamma$ acts trivially on $Ker\tau.$

\end{thm}

$Proof.$ Let $B_{\Gamma(\tau)}$ denote the normal subgroup of B generated by the elements $^{\gamma}gg^{-1}$ such that $\tau(^{\gamma}gg^{-1}) = 1, \gamma\in \Gamma, g\in G.$ It is evident that $B_{\Gamma(\tau)}$ is a $\Gamma$-subgroup of B and the canonical map $\delta:B\rightarrow B/B_{\Gamma(\tau)}$ is a $\Gamma$-homomorphism.

Let $\alpha: G\rightarrow B$ be a section map for the $\Gamma$-homomorphism $\tau$, then the sequence
$$
E_{\Gamma}: 1\rightarrow A/B_{\Gamma(\tau)}\rightarrow B/B_{\Gamma(\tau)}\overset{\tau'}\rightarrow G\rightarrow 1
$$
is a $\Gamma$-extension of $\Gamma$-groups with $\Gamma$-section map. where $\tau'$ is induced by $\tau$. In effect, for that it suffices to show that if $^{\sigma}g= ^{\gamma}g$  then $^{\sigma}(\delta\alpha(g)) = ^{\gamma}(\delta\alpha(g))$, $\sigma, \gamma \in \Gamma, g\in G$. One has $\tau(^{\sigma}\alpha(G)) = ^\sigma{(\tau\alpha(g))} = ^{\sigma}g$  and  $\tau(^{\gamma}\alpha(G)) = ^\gamma{(\tau\alpha(g))} = ^{\gamma}g$. Thus $\tau(^{\sigma}\alpha(G)) = \tau(^{\gamma}\alpha(G))$ and therefore $\tau(^{\gamma^{-1}\sigma}\alpha(g)) = ^{\gamma^{-1}}\tau(^{\sigma}\alpha(g)) = ^{\gamma^{-1}}\tau(^{\tau}\alpha(g)) = g.$ It follows that $^{\gamma^{-1}\sigma}\alpha(g) . \alpha(g)^{-1} \in B_{\Gamma(\tau)}$ implying the equality $\delta (^{\gamma^{-1}\sigma}\alpha(g)) = \delta (\alpha(g))$ and finally the required equality.

Now assume the sequence E possesses the $\Gamma$-property. Then the equali $B_{\Gamma(\tau)} = 1$ implies the isomorphism of the sequences $E \cong E_{\Gamma}$. Conversely let the sequence E satisfies the conditions of the theorem. That means it admits a $\Gamma$-section map $\alpha:G\rightarrow B$ and the group $\Gamma$ acts trivially on $Ker\tau.$ Let $\tau(^{\gamma}b\cdot b^{-1}) = 1$ for some $b\in B, \gamma\in \Gamma.$ This yields the equality $^{\gamma}(\alpha\tau(b)) = \alpha\tau(b).$ By using the equality $b = \alpha\tau(b)\cdot c, c\in A$, one obtains $^{\gamma}b = ^{\gamma}(\alpha\tau(b))\cdot ^{\gamma}c = \alpha\tau(b)\cdot c = b.$ This completes the proof.

\begin{cor}
The sequence $E_{\Gamma}$ possesses the $\Gamma$-property and every its section map is a $\Gamma$-section map.
\end{cor}
\begin{defn}
A $\Gamma$-group G is called $\Gamma$-perfect if $G = [G,G]_{\Gamma}$ or equivalently, if $H^{\Gamma}_1(G) = 0$ [33,27].
\end{defn}

{Example}

Let F(G) be the $\Gamma$-free group generated by the $\Gamma$-group G. The short exact sequence of $\Gamma$-groups
$$
1\rightarrow R\rightarrow F(G)\overset{\tau}\rightarrow G\rightarrow 1,
$$
where $\tau(\mid g \mid) = g$, has a $\Gamma$-section map $\alpha: G\rightarrow F(G)$, sending any element g to $\mid g \mid$. Then the short sequence of $\Gamma$-groups
\begin{equation}
0\rightarrow R/[F(G),R]_{\Gamma}\rightarrow F(G)/[F(G),R]_{\Gamma}\overset{\tau'}\rightarrow G\rightarrow 1
\end{equation}
 is a central $\Gamma$-equivariant extension of G having a $\Gamma$-section map $\eta\alpha$ and $\Gamma$ is trivially acting on  $R/[F(G),R]_{\Gamma}$, where $\tau'$ is induced by $\tau$ and $\eta:F(G)\rightarrow  F(G)/[F(G),R]_{\Gamma}$ is the canonical $\Gamma$-homorphism. Therefore by Theorem 3.2 the sequence (3.1) has the $\Gamma$-property.

If the $\Gamma$-group G is $\Gamma$-perfect the sequence (3.1) yields the following $\Gamma$-extension of G
\begin{equation}
  0\rightarrow R\cap[F(G),F(G)]_{\Gamma}/[F(G),R]_{\Gamma}\rightarrow [F(G),F(G)]_{\Gamma}/[F(G),R]_{\Gamma}\overset{\tau''}\rightarrow G\rightarrow 1.
 \end{equation}
 Since the group $[F(G),F(G)]_{\Gamma}/[F(G),R]_{\Gamma}$ is a $\Gamma$-subgroup of $F(G)/[F(G),R]_{\Gamma}$, the sequence (3.2) also has the $\Gamma$-property and therefore it is a $\Gamma$-extension of G (with $\Gamma$-section map). The sequence (3.2) is the universal central $\Gamma$-equivariant extension of the $\Gamma$-perfect group G and the group $R\cap[F(G),F(G)]_{\Gamma}/[F(G),R]_{\Gamma}$ is isomorphic to $H^{\Gamma}_{2}(G)$ [27]. As we see in this example the $\Gamma$-property has been used substantially.

 Now we continue our investigation of $\Gamma$-extensions of $\Gamma$-groups by considering the non-abelian case. Let
 \begin{equation}
 1\rightarrow J\overset{\sigma}\rightarrow X\overset{\tau}\rightarrow G\rightarrow 1
 \end{equation}
be an extension of $\Gamma$-groups having the $\Gamma$-property. By Theorem 3.2 the group $\Gamma$ acts trivially on J and any section map for the sequence (3.3) is a $\Gamma$-map. By conjugation one gets a $\Gamma$-homomorphism $\theta:X\rightarrow AutJ$ implying the $\Gamma$-homorphism $\psi:G\rightarrow AutJ/InJ$, where $\Gamma$ is assuming trivially acting on AutJ and InJ denotes the group of inner automorphisms of J.
\begin{defn}
The triple $(G,J,\psi)$ is called abstract kernel of the non-abelian extension (3.3) of $\Gamma$-groups.
\end{defn}
\begin{thm}
1) For any abstract kernel $(G,J,\psi)$ there is a correctly defined element,called obstruction for $(G,J,\psi)$, and belonging to $H^{3}_{\Gamma}(G,C)$. The abstract kernel $(G,J,\psi)$ possesses an extension iff $Obs(G,J,\psi) = 0$.

2) If there exists an extension of the $\Gamma$-group G with abstract kernel $(G,J,\psi)$, then the set of equivalence classes of extensions with $\Gamma$-property of G by J is bijective with $H^{2}_{\Gamma}(G,C)$, where C is the center of J.
\end{thm}
$Proof.$ We will follow the classical proof (when the action of $\Gamma$ on G is trivial).

First of all it should be noted that for a given extension (3.3) of the abstract kernel $(G,J,\psi)$ a section map $\alpha:G\rightarrow X$ with $\alpha(1) = 1$ induces by conjugation an automorphism $\varphi(x)\in \psi(x), x\in G$ of the group J and maps $f,\mu:G\times G\rightarrow J$ satisfying the well known equalities

$\alpha(x)+l = \varphi(x)(l) + \alpha(x)$,$l\in J, x\in G,$

$\alpha(x) + \alpha(y) = f(x,y) + \alpha(xy), x,y\in G,$

(1) $\varphi(x)(f(x,y)) + f(x,yz) = f(x,y) + f(xy,z), z\in G,$

(2) $\varphi(x)\varphi(y) = \mu(f(x,y))\varphi(xy).$

Taking into account the action of $\Gamma$ one has $^{\gamma}(\alpha(x)+^{\gamma}l = ^{\gamma}(\varphi(x)(l)) + ^{\gamma}(\alpha(x))$. Therefore $\alpha(^{\gamma}(x)) + l = ^{\gamma}\varphi(x)(l) + \alpha(^{\gamma}x)$. On the other hand $\alpha(^{\gamma}x)+l = \varphi(^{\gamma}x)(l) + \alpha(^{\gamma}x).$ Thus $\varphi(^{\gamma}x)(l) = ^{\gamma}\varphi(x)(l)$ showing that $\varphi$ is a $\Gamma$-map. Similarly it can be proved that the maps f and $\mu$ are $\Gamma$-maps.

Conversely, for given $\Gamma$-maps $\varphi:G\rightarrow Aut(J), f,\mu:G\times G\rightarrow J$ satisfying the equalities (1) and (2), and $\varphi(1) = 1, f(1,g) = f(x,1) = 0, x\in J, g\in G,$ we can construct an extension of the $\Gamma$-group G with $\Gamma$-property by considering the set $B(J,\varphi,f,G)$ of couples $(x,g), x\in J, g\in G,$ and defining the group structure on it as follows:

$(x,g) + (x_1,g_1) = (x + \varphi(g)x_1 + f(g,g_1),gg_1)$, and $\Gamma$ is acting on $B(J,\varphi,f,G)$ componentwise $^{\gamma}(x,g) = (x,^{\gamma}g), \gamma \in \Gamma$.This $\Gamma$-extension is called semi-direct product $\Gamma$-extension of the $\Gamma$-group G by J. It is evident that it has the $\Gamma$-property.

For any $\Gamma$-extension of G by J with $\Gamma$-property and abstract kernel $(G,J,\psi)$ the arising $\Gamma$-maps $\varphi$, f and $\psi$ should satisfy the equalities (1) and (2).That is not the case in  general and the obstruction is defined by the equality

(3) $\varphi(x)(f(x,y)) + f(x,yz) = k(x,y,z) + f(x,y) + f(xy,z),$

where k(x,y,z) is an element of the center C of J which is the kernel of $\mu.$ It is evident that $k:G\times G\times G\rightarrow C$ is a $\Gamma$-map and it can be proved similarly to the classical case that it is 3-th $\Gamma$-cocycle of the chain complex $C^{*}_{\Gamma}(G,C).$ The $\Gamma$-cocycle k is called obstruction $Ob(G,J,\psi)$ for the abstract kernel $(G,J,\psi).$

Finally, any $\Gamma$-extension of the abstract kernel $(G,J,\psi)$ is equivalent to a semi-direct product $\Gamma$-extension and the set of semi-direct product $\Gamma$-extensions of G by J is bijective to $H^{2}_{\Gamma}(G,C)$. The proof of theses two assertions completely follows the well known case when $\Gamma$ is acting trivially and it is left to the reader. This completes the proof.

As noted in Preliminaries the bijection of the set of $\Gamma$-equivariant extensions of G by A with the second $\Gamma$-equivariant homology group of the $\Gamma$-group G could be extended to higher dimensions. To this aim we need the notion of n-fold $\Gamma$-equivariant extension of G by A which is defined as a long exact sequence of $\Gamma$-groups:
  $$
  0\rightarrow A\rightarrow B_{1}\overset{\alpha_{1}}\rightarrow B_{2}\rightarrow ... \rightarrow B_{n}\overset{\alpha_{n}}\rightarrow G\rightarrow 1,
  $$
  where $0\rightarrow A\rightarrow B_{1}\rightarrow Im\alpha_{1}\rightarrow 0$, $Im\alpha_{i-1}\rightarrow   B_{i}\rightarrow Im\alpha_{i}\rightarrow\rightarrow 0$, $2\leq i\leq n-1$, are proper sequences of $\Gamma$-equivariant G-modules and $0\rightarrow Im\alpha_{n-1}\rightarrow B_{n}\rightarrow G\rightarrow e$ is a $\Gamma$-equivariant extension of G by $Im\alpha_{n-1}$.

  Let $E^{n}_{\Gamma}(G,A)$, $n\geq 1$, denote the class of equivalence classes of n-fold $\Gamma$-equivariant extensions of G by A. It becomes an abelian group by using the Baer sum operation for all $n\geq 1$.
  \begin{thm}
   There is an isomorphism
 $E^{n}_{\Gamma}(G,A)\cong  H^{n+1}_{\Gamma}(G,A)$ for $n\geq 1.$
  \end{thm}
  $Proof.$ We need the following property of the relative derived functors $Ext^{*}_{\Gamma}(A,B)$ of $Hom_{G\propto \Gamma}(A,B)$ defined in the category of $\Gamma$-equivarianr G-modules by using relatively projective $\Gamma$-resolutions with respect to the class $\mathcal{P}$ of proper sequences.

   Every proper sequence of $\Gamma$-equivariant G-modules $0\rightarrow A\rightarrow D\rightarrow B\rightarrow 0 $ and any $\Gamma$-equivariant G-module L give rise long exact sequences

$0\rightarrow Ext^{0}_{\Gamma}(B,L)\rightarrow Ext^{0}_{\Gamma}(D,L)\rightarrow Ext^{0}_{\Gamma}(A,L))\overset{\theta^{0}}\rightarrow Ext^{1}_{\Gamma}(B,L)\rightarrow ... \rightarrow Ext^{n}_{\Gamma}(A,L))\overset{\theta^{0}}\rightarrow Ext^{n+1}_{\Gamma}(B,L)\rightarrow ... ,$

$0\rightarrow Ext^{0}_{\Gamma}(L,A)\rightarrow Ext^{0}_{\Gamma}(L,D)\rightarrow Ext^{0}_{\Gamma}(L,B))\overset{\delta^{0}}\rightarrow Ext^{1}_{\Gamma}(L,A)\rightarrow ... \rightarrow Ext^{n}_{\Gamma}(L,B))\overset{\delta^{0}}\rightarrow Ext^{n+1}_{\Gamma}(L,A)\rightarrow ... .$

The connected sequence of contravariant functors $(Ext^{n}_{\Gamma}(-,A),\theta^{n},n=1,2, ...)$ is the right universal sequence of contravariant functors (or right satellite of of the functor $Ext^{1}_{\Gamma}(-,L)$) relative to the class $\mathcal{P}$ of proper sequences of $\Gamma$-equivariant G-modules [26]. In effect, let $(U_{n}, \mu^{n}, n=1,2, ...)$ be a connected sequence of contravariant functors related to the class $\mathcal{P}$ and $f^{1}:Ext^{1}_{\Gamma}(-,A)\rightarrow U^{1}(-)$ be a morphism of functors. It will be shown that there exists a unique extension of the morphism $f^{1}$ to $f^{n}:Ext^{n}_{\Gamma}(-,A)\rightarrow U^{n}(-), n=2,3, ...$ compatible with the connecting homomorphisms. To define the extension $f^{2}: Ext^{2}(J,A)\rightarrow U^{2}(J)$ consider the short exact sequence $0\rightarrow L\rightarrow P\rightarrow J\rightarrow0$, where P is a relatively projective $\Gamma$-equivariant G-module, implying the isomorphism $\theta^{1}:Ext^{1}_{\Gamma}(L,A)\rightarrow Ext^{2}_{\Gamma}(J,A)$. Then the homomorphism $f^{2}$ is given by $\mu^{1}f^{1}(\theta^{1})^{-1}$. It is easily checked that $f^{2}$ is correctly defined, compatible with the connecting homomorphisms $\theta^{1}$ and $\delta^{1}$, and it is the unique extension of $f^{1}$. The extension of $f^{n}$ to $f^{n+1}$ for $n> 1$ is constructed similarly.

Now it will be shown that the connected sequence of contravariant functors $(E^{n}_{\Gamma}(-,A),\sigma^{n},n=1,2, ...)$ in the category of $\Gamma$-equivariant G-modules is also the right universal sequence of contravariant functors relative to the class $\mathcal{P}$ of proper sequences of $\Gamma$-equivariant G-modules. For the proper sequence $E= 0\rightarrow J\rightarrow B\rightarrow L\rightarrow0$ the connecting homomorphism $\delta^{n}: E^{n}_{\Gamma}(J,A)\rightarrow E^{n+1}_{\Gamma}(L,A)$ is given by $\sigma^{n}([E_{n}]) = [E_{n} \otimes E]$, where
$E_{n} = 0\rightarrow A\rightarrow X_{1}\rightarrow X_{2}\rightarrow ... \rightarrow X_{n}\rightarrow J\rightarrow 0$ and $E_{n} \otimes E$ is an n+1- fold extension of L by A obtained by splicing the n-fold extension $E_{n}$ with the proper sequence E.

Let $\varphi^{1}:E^{1}_{\Gamma}(-,A)\rightarrow U^{1}(-)$ be a morphism of functors. Its uniquely defined extension  $\varphi^{2}:E^{1}_{\Gamma}(-,A)\rightarrow U^{1}(-)$ is realized as follows: for $[E_{2}] \in Ext^{2}_{\Gamma}(J,A)$, $E_{2} = 0\rightarrow A\rightarrow X_{1}\overset{\alpha_{1}}\rightarrow X_{2}\rightarrow J\rightarrow 0$, by using connecting homomorphisms with respect to the proper sequence $0\rightarrow Im\alpha_{1}\rightarrow X_{2}\rightarrow J\rightarrow 0$ we define $\varphi^{2}([E_{2}]) = \mu^{1}\varphi^{1}([E_{1}])$, where $E^{1} = 0\rightarrow A\rightarrow X_{1}\rightarrow Im\alpha_{1}\rightarrow 0$. The extension of $\varphi^{n}$ to $\varphi^{n+1}$ for $n> 1$ is constructed in a completely similar way and it is omitted.

We conclude that the isomorphism $f^{1}: Ext^{1}_{\Gamma}(-,A)\rightarrow E^{1}_{\Gamma}(-,A)$ implies the isomorphism of these two right universal sequences of contravariant functors and this yields the isomorphism $Ext^{n}_{\Gamma}(B,A) \cong E^{n}_{\Gamma}(B,A)$ for $n\geq 1$ and $\Gamma$-equivariant G-modules A and B. It follows that there is an isomorphism $h^{1}: E^{1}_{\Gamma}(G,-)\rightarrow E^{2}_{\Gamma}(\mathbb{Z},-)$.

If we consider the functors $E^{n}_{\Gamma}(B,A)$ with respect to the second variable then the sequence $E^{n}_{\Gamma}(\mathbb{Z},-), \eta^n, n\geq 2$, and the sequence $E^{n}_{\Gamma}(G,-),\lambda^{n}, n\geq 1$, are right universal sequences of covariant functors on the category of $\Gamma$-equivariant G-modules with respect to the class $\mathfrak{P}$ of proper sequences of $\Gamma$-equivariant G-modules. The connecting homomorphisms $\eta^{n}$ and $\lambda^{n}$ are defined and the universality is shown similarly to the case of the previous sequence $(E^{n}_{\Gamma}(-,A),\sigma^{n},n=1,2, ...)$. Therefore the isomorphism $h^{1}: E^{1}_{\Gamma}(G,-)\rightarrow E^{2}_{\Gamma}(\mathbb{Z},-)$ induces isomorphisms $h^{n}: E^{n}_{\Gamma}(G,-)\rightarrow E^{n+1}_{\Gamma}(\mathbb{Z},-)$ for all $n\geq 1$. It remains to apply the isomorphism $E^{n+1}_{\Gamma}(\mathbb{Z},-)\rightarrow Ext^{n+1}_{\Gamma}(\mathbb{Z},-)$. This completes the proof.

\;

\section{Some computations}

It is reasonable to ask for the computation of the $\Gamma$-equivariant (co)homology of groups introduced in [27]. As noted above $\Gamma$-equivariant homology and cohomology groups of retracts of $\Gamma$-free groups are trivial for $n> 1$. Here we provide an attempt to the investigation of this problem for finite cyclic $\Gamma$-groups.

It is well known that for the computation of the (co)homology of finite cyclic groups the following free resolution of $\mathds{Z}$ is used:
\begin{equation}
 ... \overset{D}\rightarrow \mathds{Z}(\mathrm{\mathds{Z}_{m}})\overset{N}\rightarrow \mathds{Z}(\mathrm{\mathds{Z}_{m}})\overset{D}\rightarrow \mathds{Z}(\mathrm{\mathds{Z}_{m}})\overset{N}\rightarrow  \mathds{Z}(\mathrm{\mathds{Z}_{m}})\overset{D}\rightarrow \mathds{Z}(\mathrm{\mathds{Z}_{m}})\overset{\epsilon}\rightarrow \mathds{Z}\rightarrow 0,
\end{equation}
where $\mathrm{\mathds{Z}_{m}}$ is a finite cyclic group of order m and generator t, D = t - 1 and $N = 1 + t + t^{2} + ... + t^{m-1}$. Assume now that a group $\Gamma$ is acting on $\mathds{Z}_{m}$ with trivial action on $\mathds{Z}$. Similarly to the action of $\Gamma$ on the bar resolution of $\mathds{Z}$ it induces an action of $\Gamma$ on the resolution (4.1) of  $\mathds{Z}$. It is easily checked that the homomorphism D is compatible with trivial action of $\Gamma$ only. This case is not interesting, since it is reduced to the usual (co)homology of groups. Therefore  the resolution (4.1) is unsuitable for the computation of $\Gamma$-equivariant (co)homology of finite cyclic $\Gamma$-groups.

By slightly changing the value of D we are able to compute the rational $\Gamma$-equivariant (co)homology of the finite cyclic $\Gamma$-group $\mathds{Z}_{m}$.

Let $B_{*} \otimes \mathds{Q}\rightarrow \mathds{Q}$ be the bar resolution of the field $\mathds{Q}$ of rational numbers obtained by tensoring the bar resolution $B_{*}\rightarrow \mathds{Z}$ of $\mathds{Z}$ by $\mathds{Q}$, where $B_{0} \otimes \mathds{Q} \cong \mathds{Q}(G)$, and $B_{n}, n\geq 1,$ is the free $\mathds{Q}(G)$- module generated by the elements $[g_{1}, g_{2}, ..., g_{n}], g_{i}\in G$.

The rational (co)homology of groups with coefficients in $\mathds{Z}(G)$-modules is defined as follows:
\begin{defn}

           $_{\mathds{Q}}H_{n}(G,A) = H_{n}((B_{*}\otimes \mathds{Q})\otimes_{\mathds{Z}(G)} A)$ and  $_{\mathds{Q}}H^{n}(G,A) =
           Hom_{\mathds{Z}(G)}(B_{*}\otimes \mathds{Q}, A)$ for $n\geq 0$ and any $\mathds{Z}(G)$-module A

\end{defn}

It is easily checked that so defined rational homology and cohomology of groups don't depend on the projective $\mathds{Q}(G)$-resolution of $\mathds{Q}$ and there are isomorphisms $_{\mathds{Q}}H_{n}(G,A) \cong H_{n}(G,A\otimes \mathds{Q})$ and $_{\mathds{Q}}H^{n}(G,A) = H^{n}(G, Hom(\mathds{Q},A))$, $n\geq 0$. In particular taking into account the isomorphisms $\mathds{Q} \otimes \mathds{Q} \cong \mathds{Q}$  and $Hom(\mathds{Q},\mathds{Q}) \cong \mathds{Q}$ one has the isomorphisms $_{\mathds{Q}}H_{n}(G,\mathds{Q}) \cong H_{n}(G,\mathds{Q})$ and $_{\mathds{Q}}H^{n}(G,\mathds{Q}) \cong H^{n}(G,\mathds{Q})$, $n\geq 0$, for the trivial $\mathds{Z}(G)$-module $\mathds{Q}$.

Based on this definition of rational (co)homology of groups we provide the following definition of rational $\Gamma$-equivariant homology and cohomology of $\Gamma$-groups. Let $\Gamma$ be a group acting on the group G and trivially acting on $\mathds{Q}$ implying the action of $\Gamma$ on the bar resolution $B_{*}\otimes \mathds{Q}$.
\begin{defn}

 $_{\mathds{Q}}H^{\Gamma}_{n}(G,A) = H_{n}((B_{*}\otimes \mathds{Q})\otimes_{\mathds{Z}(G\rtimes \Gamma)} A))$ and  $_{\mathds{Q}}H_{\Gamma}^{n}(G,A) = Hom_{\mathds{Z}(G\rtimes \Gamma)}(B_{*}\otimes \mathds{Q}, A)$ for $n\geq 0$ for any $\Gamma$-equivariant $\mathds{Z}(G\rtimes \Gamma)$-module A

\end{defn}

We recall that any $t^{i}$ is a generator of $\mathds{Z}_{m}$ iff (i,m)= 1 and the total number of generators is equal to the number of integers coprime with m and less than m which is called the Euler $\varphi$-function. Any automorphism of $\mathrm{Z}_{m}$ maps a generator to a generator and therefore it has the form $t\rightarrow t^{k}$, $1\leq k< m$, where (k,m)= 1, k and m not having non-trivial commun divisors.

 Since $\Gamma$ is not trivially acting, that means there is $\gamma \in\Gamma$ such that $^{\gamma}t = t^{k}$, $k> 1$. For the computation of rational (co)homology of $\mathrm{Z}_{m} = \{1,t,t^{2},...,t^{m-1}\} $, $t^{m} = 1$, we will use the following sequence of free $\mathrm{Q}(\mathrm{Z}_{m})$-modules over $\mathrm{Q}$:
\begin{equation}
 ... \overset{D}\rightarrow \mathds{Q}(\mathds{Z}_{m})\overset{N}\rightarrow \mathds{Q}(\mathds{Z}_{m})\overset{D}\rightarrow \mathds{Q}(\mathds{Z}_{m})\overset{N}\rightarrow  \mathds{Q}(\mathds{Z}_{m})\overset{D}\rightarrow \mathds{Q}(\mathds{Z}_{m})\overset{\epsilon}\rightarrow \mathds{Q}\rightarrow 0,
\end{equation}
where $D = t^{m-1}+t^{m-2}+...+t - (m-1) $ and $N = 1+t+t^{2}+ ...+ t^{m-1}$. It will be shown that the sequence (4.2) is a projective $\mathrm{Q}(\mathrm{Z}_{m})$-resolution of $\mathds{Q}$. It is evident that the homomorphims D and N are compatible with the action of $\Gamma$ and DN = 0, ND = 0, $\epsilon$D  = 0.

 For $f(t) = q_{0} + q_{1}t + q_{2}t^{2}+ ... + q_{m-1}t^{m-1} \in \mathds{Q}(\mathds{Z}_{m})$, $q_i\in Q$, one has

 $ D(f(t)) = -(m-1)q_{0} + q_{1} + q_{2} + ... + q_{m-1} + (q_{0}-(m-1)q_{1} + q_{2} + ... + q_{m-1})t + ... + (q_{0} + q_{1} + ... + q_{m-3}  -(m-1)q_{m-2}+ q_{m-1})t^{m-2} + (q_{0} + q_{1} + ... + q_{m-2}-(m-1) q_{m-1})t^{m-1}$.

If $D(f(t)) = 0$ this yields a system of m equalities: $q_{0} + q_{1} + ... + q_{i-1} -(m-1)q_{i} + q_{i+1}+ .. + q_{m-1} = 0, i = 0,1,...,m-1$, implying the equalities $q_{0} = q_{1} = q_{2} = ... = q_{m-2} = q_{m-1}$ and therefore $N(q_{0}) = f(t)$.

If N(f(t)) = 0, then $q_{0} + q_{1} + ... + q_{m-2} + q_{m-1} = 0$. Assume there is $\varphi(t) = x_{0} + x_{1}t + x_{2}t^{2} + ... + x_{m-1}t^{m-1}$ such that $D(\varphi(t)) = f(t)$. This yields the equalities

$-(m-1)x_{0}+ x_{1} + x_{2} + ... + x_{m-1} = q_{0}, $

$x_{0} - (m-1)x_{1} + x_{2}  + ... + x_{m-1} = q_{1},$

$....................................................$

$x_{0} + x_{1} + x_{2}  + ... + x_{m-2} -(m-1)x_{m-1} = q_{m-1}$.

Since $\Sigma_{i}q_{i} = 0$, this system of m linear equations has infinitely many solutions in the field $\mathrm{Q}$ of rational numbers, and for every $q\in Q$ by taking $x_m = q$ the solution has the form $x_{i} = (q_{m-1} - q_{i})/m + q, i = 0,1,...,m-2 $, $x_{m-1} = q$. The same holds for the case $\epsilon(f(t)) = 0$. It follows that the sequence (4.2) is exact.

It remains to show that $0\rightarrow KerD\rightarrow \mathds{Q}(\mathds{Z}_{m})\rightarrow ImD\rightarrow 0$ and $0\rightarrow KerN\rightarrow \mathrm{Q}(\mathrm{Z}_{m})\rightarrow ImN\rightarrow 0$ are sequences of $\Gamma$-equivariant $\mathds{Q}(\mathds{Z}_{m})$-modules having $\Gamma$-section map. Since ImN = KerD, any element of ImN has the form $f(t) = q + qt + qt^{2} + ... + qt^{m-1}$ and $\Gamma$ acts trivially on f(t). Therefore the map $\alpha: ImN\rightarrow \mathrm{Q}(\mathrm{Z}_{m})$, $\alpha(f(t)) = q$ is a $\Gamma$-section map for the sequence $0\rightarrow KerN\rightarrow \mathds{Q}(\mathds{Z}_{m})\rightarrow ImD\rightarrow 0$. For the case of the sequence $0\rightarrow KerD\rightarrow \mathds{Q}(\mathds{Z}_{m})\rightarrow ImD\rightarrow 0$ every element $f(t)\in ImD$ satisfies the equality $ \Sigma_{i}q_{i} = 0$ since ImD = KerN, and consider the element $\varphi(t) = -1/m (f(t))$. It will be shown that $D(\varphi(t)) = f(t)$. In effect one has $D(\varphi(t)) =$

$(-(m-1)(-1/m (q_{0}))+ q_{1} + q_{2} + ... + q_{m-1}) + $

$(-1/m(q_{0}) - (m-1)(-1/m(q_{1})) -1/m (q_{2})  - ... -1/m (q_{m-1}))t + $

$+ ... + (-1/m(q_{0})  -1/m(q_{1}) - 1/m(q_{2})  - ... -1/m( q_{m-2}) -(m-1)(-1/m(q_{m-1})))t^{m-1} = $

$q_{0} -1/m(q_{0} + q_{1} + q_{2} + ... + q_{m-1}) + $

$(q_{1} - (-1/m(q_{0} + q_{1} + q_{2} + ... + q_{m-1}))t + $

$+ ... + (q_{m-1} - (-1/m(q_{0} + q_{1} + q_{2}  + ... + q_{m-2}+ q_{m-1})t^{m-1} = f(t)$.

Therefore the map $\alpha: ImD\rightarrow \mathrm{Q}(\mathrm{Z}_{m})$ sending f(t) to -1/m(f(t)) is a $\Gamma$-section map.

  We have proven that the sequence (4.2) is a $\mathds{Q}(\mathds{Z}_{m})$-equivariant projective resolution of the trivial $\mathds{Q}(\mathds{Z}_{m})$- module $\mathds{Q}$, where $\mathds{Q}(\mathds{Z}_{m})$ is a relatively free  $\mathds{Q}(\Gamma)$-equivariant $\mathds{Z}(\mathds{Z}_{m}\rtimes \Gamma)$-module. Therefore the homology groups of the complex

$...\overset{D_{*}}\rightarrow\mathds{Q}(\mathds{Z}_{m})\otimes_{\mathds{Z}(\mathrm{Z}_{m}\rtimes \Gamma)} A \overset{N_{*}}\rightarrow  \mathds{Q}(\mathds{Z}_{m})\otimes_{\mathds{Z}(\mathds{Z}_{m}\rtimes \Gamma)} A \overset{D_{*}}\rightarrow \mathds{Q}(\mathds{Z}_{m})\otimes_{\mathds{Z}(\mathds{Z}_{m}\rtimes \Gamma)} A \overset{N_{*}}\rightarrow  \mathds{Q}(\mathds{Z}_{m})\otimes_{\mathds{Z}(\mathds{Z}_{m}\rtimes \Gamma)} A \overset{D_{*}}\rightarrow \mathds{Q}(\mathds{Z}_{m})\otimes_{\mathds{Z}(\mathds{Z}_{m}\rtimes \Gamma)} A\rightarrow 0$

give us the rational $\Gamma$-equivariant homology $_{\mathds{Q}}H_{n}(\mathrm{Z}_{m}, A)$ of the cyclic group $\mathds{Z}_{m}$ with coefficients in the $\mathds{Z}((\mathrm{Z}_{m}\rtimes \Gamma)$-module A. By applying the isomorphisms $ \mathds{Q}(\mathds{Z}_{m})\otimes_{\mathds{Z}(\mathds{Z}_{m}\rtimes \Gamma)} A \cong (\mathds{Q}\otimes \mathds{Z}(\mathds{Z}_{m}))\otimes_{\mathds{Z}(\mathrm{Z}_{m}\rtimes \Gamma)} A \cong \mathds{Q}\otimes (\mathds{Z}(\mathds{Z}_{m})\otimes_{\mathds{Z}(\mathds{Z}_{m}\rtimes \Gamma)} A) \cong \mathds{Q}\otimes A_{\Gamma}$ and the fact that every element of $\mathds{Q}\otimes A_{\Gamma}$ can be written in the form $\sum_i(q_i\otimes [a_i])$, $q_i\in \mathds{Q}$, $[a_i]\in A_{\Gamma}$. we have finally obtained

\begin{thm}

Let $\Gamma$ be a group not trivially acting on the finite cyclic group $\mathds{Z}_{m}$. Then for any $\Gamma$-equivariant $\mathds{Z}(\mathds{Z}_{m}\rtimes \Gamma)$-module A

$_{\mathds{Q}}H^{\Gamma}_{0}(\mathds{Z}_{m},A) = \mathds{Q}\otimes A_{\Gamma} $

$_{\mathds{Q}}H^{\Gamma}_{2n-1}(\mathds{Z}_{m},A) =  \sum_i(q_i\otimes [a_it^{m-1} + a_it^{m-2} + ... + a_it]) = \sum_i(q_i\otimes [(m-1)a_i])/ Im N_{*} $

$_{\mathds{Q}}H^{\Gamma}_{2n}(\mathds{Z}_{m},A) = \sum_i(q_i\otimes [a_it^{m-1} + a_it^{m-2} + ... + a_it + a_i]) = 0)/ Im D_{*} $,

for $n> 0$ and the homomorphisms $D_{*}$ and $N_{*}$ are induced by D and N respectively.

\end{thm}

For the rational $\Gamma$-equivariant cohomology $_{\mathds{Q}}H^{n}(\mathrm{Z}^{m}, A)$ of the cyclic group $\mathds{Z}_{m}$ with coefficients in the $\mathds{Z}((\mathds{Z}_{m}\rtimes \Gamma)$-module A we consider the complex

$0\rightarrow Hom_{\mathds{Z}(\mathrm{Z}_{m}\rtimes \Gamma)}(\mathds{Q}(\mathrm{Z}_{m}), A)\overset{D^{*}}\rightarrow Hom_{\mathds{Z}(\mathrm{Z}_{m}\rtimes \Gamma)}(\mathds{Q}(\mathds{Z}_{m}), A)\overset{N^{*}}\rightarrow Hom_{\mathds{Z}(\mathrm{Z}_{m}\rtimes \Gamma)}(\mathds{Q}(\mathds{Z}_{m}), A)$

$\overset{D^{*}}\rightarrow Hom_{\mathds{Z}(\mathrm{Z}_{m}\rtimes \Gamma)}(\mathds{Q}(\mathds{Z}_{m}), A) \overset{N^{*}}\rightarrow Hom_{\mathds{Z}(\mathds{Z}_{m}\rtimes \Gamma)}(\mathds{Q}(\mathds{Z}_{m}), A)\overset{D^{*}}\rightarrow ... $

\begin{thm}

Let $\Gamma$ be a group not trivially acting on the finite cyclic group $\mathrm{Z}_{m}$. Then for any $\Gamma$-equivariant $\mathds{Z}(\mathds{Z}_{m}\rtimes \Gamma)$-module A

$_{\mathds{Q}}H_{\Gamma}^{0}(\mathds{Z}_{m},A) = Hom(\mathds{Q},A^{\Gamma}) $

$_{\mathds{Q}}H_{\Gamma}^{2n-1}(\mathds{Z}_{m},A) = Ker N^{*} / Im D^{*} $

$_{\mathds{Q}}H_{\Gamma}^{2n}(\mathds{Z}_{m},A) = Ker D^{*} / Im N^{*} $,

for $n> 0$ and the homomorphisms $D^{*}$ and $N^{*}$ are induced by D and N respectively.

\end{thm}

We see that the rational $\Gamma$-equivariant (co)homology of finite cyclic $\Gamma$-groups is of period 2 for $n> 0$.

\;

\section{$\Gamma$-derived functors and $\Gamma$-equivariant Hochschild homology}

Let $\Gamma$ be a group acting on the ring $\Lambda$ with unit, that means  a group homomorphism $\theta :\Gamma \rightarrow Aut(\Lambda)$ is given, the element $\theta (\gamma)(\lambda)$ is denoted $^{\gamma}\lambda$,  and let A be a left $\Lambda$-module on which $\Gamma$ is acting such that

$^{\gamma}(\lambda a)  =  ^{\gamma}\lambda^{\gamma}a , \gamma\in \Gamma, \lambda\in \Lambda, a\in A.$

Then A will be called $\Gamma$-equivariant left $\Lambda$-module. Denote by $\Gamma A$ the $\Lambda$-submodule of A generated by the elements ${^{\gamma}a - a}$, $a\in A$, $\gamma\in \Gamma$, and by $A_{\Gamma}$ the quotient of A by $\Gamma A$.

\begin{defn}

It will be said that a group $\Gamma$ is acting on a chain complex $L_{*}$ of left $\Lambda$-modules

$... \rightarrow L_{n}\overset{\delta_{n}}\rightarrow L_{n-1}\overset{\delta_{n-1}}\rightarrow ...$

if $\Gamma$ acts on each $L_{n}$ becoming $\Gamma$-equivariant $\Lambda$-module and every $\delta_{n}$ satisfies the following condition

(1)   $\delta_{n}(^{\gamma}l_{n}-l_{n})\in \Gamma L_{n-1}, l_{n}\in L_{n}, \gamma\in \Gamma.$

\end{defn}
In particular condition (1) is satisfied if $\delta_{n}$ is compatible with the action of $\Gamma$.

\begin{defn}

The homology groups $H^{\Gamma}_{n}(L_{*})$, $n\in \mathds{Z}$, of the chain complex $L_{*}$ are defined as the homology groups of the quotient chain complex of $L_{*}$:

$L^{\Gamma}_{*} = ... \rightarrow (L_{n})_{\Gamma}\overset{\delta'_{n}}\rightarrow (L_{n-1})_{\Gamma}\overset{\delta'_{n-1}}\rightarrow ...$

The groups $H^{\Gamma}_{n}(L_{*})$ are called $\Gamma$-equivariant homology groups of $L_{*}$.

\end{defn}

The consideration of the chain complex $L^{\Gamma}_{*}$ is motivated by the following important cases.

Case 1 -  $Relation$ $with$ $cyclic$ $homology.$

It will be said that a group $\Gamma$ is acting on a left $\Lambda$-module M if it is acting on $\Lambda$ and M such that $^{\gamma}(\lambda m) = ^{\gamma}\lambda^{\gamma}m$. In that case M is called $\Gamma$-equivariant $\Lambda$-module.

It will be said that a group $\Gamma$ is acting on unital $\kappa$-algebra A if it is acting on A and $\kappa$ such that $^{\gamma}(ka) = ^{\gamma}k^{\gamma}a$, $k\in \kappa, a\in A, \gamma \in \Gamma$. The group $\Gamma$ acts on A-bimodule M if it is acting on the $\kappa$-algebra A and on M such that $^{\gamma}((am)a') - (^{\gamma}a^{\gamma}m)^{\gamma}a'  =  ^{\gamma}a(^{\gamma}m^{\gamma}a')$. In that case M is called $\Gamma$-equivariant A-bimodule or equivalently $\Gamma$-equivariant $A^{e}$-module, where $A^{e}$ is the enveloping algebra of A, $A^{e}  =  A\otimes A^{op}$, $A^{op}$ being the opposite algebra of A.

Let
\begin{equation}
C_{*}(A,M)= ...\rightarrow M\otimes A^{\otimes n}\overset{b}\rightarrow M\otimes A^{\otimes n-1}\overset{b}\rightarrow ... \rightarrow M\otimes A \overset{b}\rightarrow M
\end{equation}
be the Hochschild complex, where the $\kappa$-module $M\otimes A^{\otimes n}$ is in degree n and the tensor product is taken over $\kappa$.

\begin{defn}
Let $\Gamma$ be a group acting on the Hochschild complex $C_{*}(A,M)$. Then  $H^{\Gamma}_{*}(C_{*}(A,M))$ is called $\Gamma$-equivariant Hochschild homology of the $\kappa$-algebra A with coefficients in $\Gamma$-equivariant $A^{e}$- module M. If the action of $\Gamma$ is induced by its actions on A and M $(^{\gamma}(m,a_{1},...,a_{n}) = (^{\gamma}m,^{\gamma} a_{1},...,^{\gamma}a_{n}))$, then $H^{\Gamma}_{*}(C_{*}(A,M))$ will be denoted $H^{\Gamma}_{*}(A,M)$ and $HH^{\Gamma}_{*}(A)$ for M = A.
\end{defn}

Let M = A and assume the group $\mathds{Z}$ of integers acts trivially on $\kappa$ and A. Define the action of $\mathds{Z}$ on $A^{\otimes n+1}$ via the composition of the canonical homomorphism $\mathds{Z}\rightarrow \mathds{Z}/(n+1)\mathds{Z} \cong \mathds{Z}_{n+1}$ with the action of $\mathds{Z}_{n+1}$ on $A^{\otimes n+1}$ given by

$t^{n}(a_{0},...,a_{n}) = (-1)^{n}(a_{n},a_{0},...,a_{n-1})$ on the generator $(a_{0},...,a_{n})$ of $A^{\otimes n+1}$, $n\geq 1$. Then $A^{\otimes n+1}$ becomes $\mathds{Z}$-equivariant $\kappa$-module and it is well known that the $\kappa$-homomorphism b satisfies condition (1) of Definition 5.1. Therefore $\mathds{Z}$ is acting on the chain complex $C_{*}(A,A)$, the chain complex $C_{*}(A,A)_{\mathds{Z}}$ is just the Connes complex and the homology groups $H^{\mathds{Z}}_{n}(C_{*}(A,A))$, $n\succeq 0$, are Connes homology groups of the $\kappa$-algebra A. It is well known that they are isomorphic to cyclic homology groups of A when $\kappa$ contains the group $\mathds{Q}$ of rational numbers.

We conclude that the cyclic homology of the algebra A over $\kappa$ containing $\mathds{Q}$ is $\mathds{Z}$-equivariant Hochschild homology $H^{\mathds{Z}}_{*}(C_{*}(A,A))$ .

Case 2 -  $Relation$ $with$ $\Gamma$-$equivariant$ $homology$ of $groups$

Let G be a group on which the group $\Gamma$ is acting and consider the bar $\mathds{Z}(G)$-resolution of $\mathds{Z}$:

$B_{*}(\mathds{Z}) = ...\rightarrow B_{n}\rightarrow B_{n-1}\rightarrow ...\rightarrow B_{1}\rightarrow B_{0}\rightarrow \mathds{Z}\rightarrow 0$,

where $B_{0} = \mathds{Z}(G)$ and $B_{n}$ is the free $\mathds{Z}(G)$ - module generated by $[g_{1},...,g_{n}]$, $g_{i}\in G, n\geq 1$. The action of the group $\Gamma$ on $B_{*}(\mathds{Z})$ is given by

$^{\gamma}(g[g_{1},...,g_{n}]) = ^{\gamma}g[^{\gamma}g_{1},...,^{\gamma}g_{n}], n\geq 1$,

and we assume $\Gamma$ is acting trivially on $\mathds{Z}$.

The action of $\Gamma$ can be extended to the integral homology complex of G:

$C_{*}(G) = ...\rightarrow C_{n}\rightarrow C_{n-1}\rightarrow ...\rightarrow C_{1}\rightarrow C_{0}\rightarrow 0$,

where $C_{0} = \mathds{Z}$, $C_{n}$ is the free abelian group generated by $[g_{1},...,g_{n}]$, $g_{i}\in G, n\geq 1$ and $^{\gamma}([g_{1},...,g_{n}]) = [^{\gamma}g_{1},...,^{\gamma}g_{n}], n\geq 1$.

According to [27] the n-th $\Gamma$-equivariant integral homology $H^{\Gamma}_{n}(G)$ of G is the n-th homology group of the chain complex $B_{*}(\mathds{Z})\otimes_{\mathds{Z}(G\rtimes \Gamma)}\mathds{Z}$. It is easily checked that $H^{\Gamma}_{n}(G)$ is isomorphic to the n-th $\Gamma$-equivariant homology $H^{\Gamma}_{n}(C_{*}(G))$ of the chain complex $C_{*}(G)$. In effect, Let $F = \Sigma_{i}\mathds{Z}(G)a_{i}$ be a $\mathds{Z}(G\rtimes \Gamma)$- module free as $\mathds{Z}(G)$-module and the set of generators $Y = \{a_{i}\}$ be a $\Gamma$-set. As mentioned above it is called relatively free $\mathds{Z}(G\rtimes \Gamma)$- module. Then the following isomorphism holds: $F\otimes_{\mathds{Z}(G\rtimes \Gamma)}\mathds{Z} \cong (\Sigma_{i}\mathds{Z} a_{i})_{\Gamma}$,implying the needed isomorphism.

Besides the group action on the Hochschild complex $C_{*}(A,A)$ of the case 1 it is interesting to consider the action of the group $\Gamma$ on $C_{*}(A,A)$ induced by the action of $\Gamma$ on the $\kappa$-algebra A, $^{\gamma}(a_{1},...,a_{n}) = (^{\gamma}a_{1},...,^{\gamma}a_{n})$ on the generators $(a_{1},...,a_{n})$ of $A^{\otimes n}, n\geq 1$. Under this action of $\Gamma$ the $\kappa$-homomorphism b is compatible and the $\kappa$-module $A^{\otimes n}$ is $\Gamma$-equivariant. We particularly mean the case when $A = \mathds{Z}(G)$ and the action of $\Gamma$ on A is induced by the action of $\Gamma$ on the group G. It is clear in this case one has isomorphisms $H^{\Gamma}_{n}(C_{*}(\mathds{Z}(G),\mathds{Z}(G))) \cong H^{\Gamma}_{n}(C_{*}(G)) \cong H^{\Gamma}_{n}(G), n\geq 0$.

Therefore the $\Gamma$-equivariant Hochschild homology contains as particular cases the cyclic homology of $\kappa$-algebras for $\mathds{Q} \subset \kappa$ and the $\Gamma$-equivariant integral homology of groups.

In what follows it will always be assumed that the $\Gamma$-equivariant Hochschild homology of the $\kappa$-algebra A is defined by the action of $\Gamma$ on A. In order to describe the $\Gamma$-equivariant Hochschild homology in terms of derived functors the notion of $\Gamma$-equivariant derived functors will be introduced.

Let $\mathds{A}^{\Gamma}_{\Lambda}$ be the category of $\Gamma$-equivariant left $\Lambda$-modules. A morphism of the category $\mathds{A}^{\Gamma}_{\Lambda}$ is a $\Lambda$-homomorphism $f:M\rightarrow M'$ such that $f(^{\gamma}m) = ^{\gamma}f(m), m\in M, \gamma\in \Gamma$. As mentioned in Preliminaries, if $\Lambda = \mathds{Z}(G)$ the category $\mathds{A}^{\Gamma}_{\Lambda}$ is equivalent to the category of $G\rtimes \Gamma$-modules. It is evident if $\Lambda = \mathrm{Z}$ and $\Gamma$ acts trivially on $\mathds{Z}$ the category $\mathds{A}^{\Gamma}_{\Lambda}$ is equivalent to the category of $\mathds{Z}(\Gamma)$-modules. A $\Gamma$-equivariant $\Lambda$-module free as $\Lambda$-module with basis a $\Gamma$-set is called relatively free $\Gamma$-equivariant $\Lambda$-module. A retract of a relatively free $\Gamma$-equivariant $\Lambda$-module is called relatively projective $\Gamma$-equivariant $\Lambda$-module. Any short exact sequence of $\Gamma$-equivariant $\Lambda$-modules having a $\Gamma$-section map is called proper exact sequence of $\Gamma$-equivariant $\Lambda$-modules.

A long exact sequence of $\Gamma$-equivariant $\Lambda$-modules

$P_{*}(M) = ...\rightarrow P_{n}\overset{\alpha_{n}}\rightarrow P_{n-1}\rightarrow ...\rightarrow P_{1}\overset{\alpha_{1}}\rightarrow P_{0}\overset{\tau}\rightarrow M\rightarrow 0$,

is called $\Gamma$-projective resolution of M, where every $P_{n}$ is relatively projective $\Gamma$-equivariant $\Lambda$-module and sequences $0\rightarrow Ker\tau\rightarrow P_{0}\overset{\tau}\rightarrow M\rightarrow0$, $0\rightarrow Ker \alpha_{n}\rightarrow P_{n}\rightarrow Im \alpha_{n}\rightarrow 0, n\geq 1$, are proper sequences.
It is obvious there is a natural action of $\Gamma$ on the chain complex $P_{*}(M)$.
\begin{defn}
  Let T be an additive covariant functor from $\mathds{A}^{\Gamma}_{\Lambda}$ to the category $\mathds{A}^{\Gamma}_{\mathds{Z}}$. The left $\Gamma$-derived functors $L^{\Gamma}_{n}T, n\geq 0$, of T are defined as $L^{\Gamma}_{n}T(M) = H^{\Gamma}_{n}(TP_{*}(M))$.
 \end{defn}

It is easily checked that these derived functors are correctly defined and they don't depend of the $\Gamma$-projetive resolution of M.

Consider the action of $\Gamma$ on the tensor product $M\otimes_{\Lambda} L$ of $\Gamma$-equivariant $\Lambda$-modules M and L induced by the action of $\Gamma$ on the couples (m,l),$^{\gamma}(m,l) = (^{\gamma}m,^{\gamma}l)$, $m\in M, l\in L$.Then $M\otimes_{\Lambda} L$ becomes a $\Gamma$-equivariant abelian group or equivalently $\mathds{Z}(\Gamma)$-module. The left $\Gamma$- derived functors of the functor $-\otimes_{\Lambda} L$: $\mathds{A}^{\Gamma}_{\Lambda}\rightarrow \mathds{A}^{\Gamma}_{\mathds{Z}}$ will be denoted $Tor^{\Lambda,\Gamma}_{n}(-,L), n\geq 0$. If $\Lambda = \mathds{Z}(G)$, we recover the functors $Tor_n^{\mathcal{P}}$ defined in [27], where $\mathcal{P}$ is the projective class of proper sequences of $\Gamma$-equivariant $\mathds{Z}(G)$-modules and L is a trivial $\Gamma$-equivariant $\mathds{Z}(G)$-module, in particular $Tor^{\mathds{Z}(G),\Gamma}_{n}(\mathds{Z},L) \cong H^{\Gamma}_{n}(G,L)$, $n\geq 0$, if $\Gamma$ acts on L trivially.

\begin{defn}
If the group $\Gamma$ is acting on the $\kappa$-algebra A, let $[A,A]_{\Gamma}$ denote the $\kappa$-submodule of A generated by the elements $\{^{\gamma}a - a, aa'- a'a\}, a,a'\in A, \gamma\in \Gamma$. It will be called the $\Gamma$-additive commutator of A.

Let A be a commutative unital $\kappa$-algebra and $\Omega^{1}(A)$ be the A-module of K$\ddot{a}$hler differentials generated by the $\kappa$-linear symbols da, $a\in A$ with defining relation $d(ab) = adb + bda, a,b\in A$. We define the action of $\Gamma$ on $\Omega^{1}(A)$ as follows: $^{\gamma}(da) = d(^{\gamma}a)), a\in A, \gamma\in \Gamma$.

The $\Gamma$-equivariant A-module $\Omega^{1}_{\Gamma}(A)$ of  K$\ddot{a}$hler differentials is defined as $(\Omega^{1}(A))_{\Gamma}$.
\end{defn}
\begin{thm}
Let $\Gamma$ be a group acting on unital $\kappa$-algebra A and on the bimodule M over A. Then one has

1. $HH^{\Gamma}_{0}(A) = A/[A,A]_{\Gamma},$

2. If A is commutative, $HH^{\Gamma}_{1}(A) \cong \Omega^{1}_{\Gamma}(A),$

3. If A is relatively projective $\Gamma$-equivariant $\kappa$-module, $H^{\Gamma}_{n}(A,M) =$ $ Tor^{A^{e},\Gamma}_{n}(A,M)$  for every $\Gamma$-equivariant A-bimodule M,

4. Morita equivalence for $\Gamma$-equivariant Hochschild homology. The inclusion maps $A\rightarrow \mathrm{M}_{r}(A), M\rightarrow \mathrm{M}_{r}(M)$, induce the isomorphism $H^{\Gamma}_{n}(A,M)\rightarrow H^{\Gamma}_{n}(\mathrm{M}_{r}(A),\mathrm{M}_{r}(M)), r\geq 1, n\geq 0$,
where $\mathrm{M}_{r}(A)$ and $\mathrm{M}_{r}(M)$ are the $\kappa$-algebra of $r\times r$-matrices over A and the module of $r\times r$-matrices over M respectively.

\end{thm}

$Proof$ 1) Taking into account the homomorphism b is a $\Gamma$-homomorphism the first equality is straightforward. If A is commutative, $HH^{\Gamma}_{0}(A) = A_{\Gamma}$.

2) It is well known that the maps $HH_{1}(A)\rightarrow \Omega^{1}(A)$ and $\Omega^{1}(A)\rightarrow HH_{1}(A)$ sending respectively the class of $a\otimes a'$ to ada' and ada' to the class of $a\otimes a'$ are inverse to each other. It suffices to remark that $^{\gamma}(a\otimes a') - a\otimes a' = ^{\gamma}a\otimes ^{\gamma}a' - a\otimes a'$ is sending to $^{\gamma}ad^{\gamma}a'$ - ada' = $^{\gamma}(ada')$ - ada' and conversely, $^{\gamma}(ada')$ - ada' is sending to $^{\gamma}(a\otimes a') - a\otimes a'$.

3) Consider the Hochschild bar complex

$C^{bar}_{*}(A) = ...\rightarrow  A^{\otimes n+1}\overset{b'}\rightarrow  A^{\otimes n}\overset{b'}\rightarrow ... \overset{b'}\rightarrow A^{\otimes 2}$,

where $A^{\otimes n+1}$ is in degree n-1, $n\geq 1$, and $b' = \sum^{n-1}_{i=0}(-1)d_{i}$, $d_{i}$ are differentials of the Hochschild complex $C_{*}(A,A)$ and the group $\Gamma$ is acting on the Hochschild bar complex as it is defined for the Hochschild complex. It is evident every $A^{\otimes n}, n\geq 2$, is relatively projective $\Gamma$-equivariant $\kappa$-module, since A possesses this property. Moreover $A^{\otimes n+2}, n\geq 0$, is relatively projective $\Gamma$-equivariant left $A^{e}$-module with the action

$(\alpha,\beta)(a_{0},...,a_{n+1}) = (\alpha a_{0},a_{1},...,a_{n},a_{n+1}\beta)$. The contracting homotopy $s: A^{\otimes n}\rightarrow A^{\otimes n+1}$, $s(a_{0},...,a_{n}) = (1,a_{0},...,a_{n})$, is a $\Gamma$-map satisfying the equality b's + sb' = id.

Therefore the chain complex $C^{bar+}_{*}(A) = C^{bar}_{*}(A)\overset{b}\rightarrow A$, $b(a\otimes a') = aa', a,a'\in A$, is a $\Gamma$-projective resolution of the $\Gamma$-equivariant $A^{e}$-module A. Upon tensoring this $\Gamma$-projective resolution with a $\Gamma$-equivariant $A^{e}$-module M one obtains the Hochschild complex because of the isomorphism $M\otimes_{A^{e}}A^{\otimes n+2} \cong M\otimes_{\kappa}A^{\otimes n}$. This implies the equalities $H^{\Gamma}_{n}(A,M) = H^{\Gamma}_{n}(C^{bar+}_{*}(A)\otimes_{A^{e}}M) = Tor^{A^{e},\Gamma}_{n}(A,M)$.

4) The action of the group $\Gamma$ on A and M induces its action on $\mathrm{M}_{r}(A)$ and $\mathrm{M}_{r}(M)$ given by$^{\gamma}[a_{ij}] = [^{\gamma}a_{ij}]$, $^{\gamma}[m_{ij}] = [^{\gamma}m_{ij}], \gamma\in \Gamma, a_{ij}\in A, m_{ij}\in M$, . This action is compatible with the natural inclusions $\mathrm{M_{r}}(A)\rightarrow \mathrm{M_{r+1}}(A)$ and $\mathrm{M}_{r}(M)\rightarrow \mathrm{M}_{r+1}(M)$ inducing the action of $\Gamma$ on $\mathrm{M}(A) = lim_{\rightarrow_{r}}\mathrm{M}_{r}(A)$ and on $\mathrm{M}(M) = lim_{\rightarrow_{r}}\mathrm{M}_{r}(M)$ respectively.

It also induces an action of $\Gamma$ on the trace map tr: $\mathrm{M}_{r}(M)\rightarrow M$, $tr([m_{ij}]) = \Sigma^{r}_{i}m_{ii}$, given by $(^{\gamma}tr)([m_{ij}]) = \Sigma{^{r}_{i}} (^{\gamma}m_{ii}), \gamma \in \Gamma$. The trace map is extended to tr: $\mathrm{M}_{r}(M)\otimes \mathrm{M}_{r}(A^{\otimes n})\rightarrow M\otimes A$. By identifying $\mathrm{M}_{r}(M)$ with $\mathrm{M}_{r}(\kappa)\otimes M$ any element of $\mathrm{M}_{r}(M)$ is a sum of elements $u_{i}v_{i}$ with $v_{i}\in \mathrm{M}_{r}(\kappa)$ and $u_{i}\in M$, and the trace map takes the form

tr $(v_{0}a_{0}\otimes ... \otimes v_{n}a_{n})$ = tr $(v_{0}...v_{n})a_{0}\otimes ... \otimes a_{n}$, $v_{i}\in \mathrm{M}_{n}(\kappa), a_{0}\in M$ and $a_{j}\in A, j\geq 1$. The action of $\Gamma$ on the extended trace map is given by $(^{\gamma}tr)(v_{0}a_{0}\otimes ... \otimes v_{n}a_{n})$ = tr $(^{\gamma}v_{0}...^{\gamma}v_{n})^{\gamma}a_{0}\otimes ... \otimes ^{\gamma}a_{n}$, $\gamma \in \Gamma$. It is evident that the extended trace map is a $\Gamma$-map taking into account that tr $(^{\gamma}(v_{0}...v_{n}))$ = tr $(^{\gamma}v_{0}...^{\gamma}v_{n})$.

Thus the extended trace map yields a morphism of chain complexes $tr_{*}: C_{*}(\mathrm{M_{r}(A), \mathrm{M_{r}(M)}})\rightarrow C_{*}(A,M)$ compatible with the action of the group $\Gamma$ and therefore a morphism $tr^{\Gamma}_{*}: C^{\Gamma}_{*}(\mathrm{M_{r}(A), \mathrm{M_{r}(M)}})\rightarrow C^{\Gamma}_{*}(A,M)$. On the other hand there is a morphism $inc^{\Gamma}_{*}: C^{\Gamma}_{*}(A,M)\rightarrow C^{\Gamma}_{*}(\mathrm{M_{r}(A), \mathrm{M_{r}(M)}})$ induced by the inclusion maps $A\rightarrow \mathrm{M}_{r}(A)$, $M\rightarrow \mathrm{M}_{r}(M)$. It is immediate that $tr^{\Gamma}_{*} inc^{\Gamma}_{*} = id$. It is well known that $inc^{\Gamma}_{*}tr^{\Gamma}_{*}$ is homotopic to id and the homotopy $h = \Sigma_{i}(-1)^{i}h_{i}$ is defined by the formula

$h_{i}(a^{0},...,a^{n}) = \Sigma e_{ij}(a^{0}_{jk})\otimes e_{11}(a^{1}_{km})\otimes...\otimes e_{11}(a^{i}_{pq})\otimes e_{1q}(1)\otimes a^{i+1}\otimes ... \otimes a^{n}$,

where the sum is extended over all possible sets of indices (j,k,m,...,p,q), $a^{0}$ is in $\mathrm{M_{r}(M)}$, other $a^{i}$ are in $\mathrm{M_{r}(A)}$ and the $e_{ij}$ denoted elementary matrices. According to this homotopy formula we have the equalities

$^{\gamma}(h_{i}(a^{0},...,a^{n})) = \Sigma ^{\gamma}e_{ij}(a^{0}_{jk})\otimes ^{\gamma}e_{11}(a^{1}_{km})\otimes...\otimes ^{\gamma}e_{11}(a^{i}_{pq})\otimes ^{\gamma}e_{1q}(1)\otimes ^{\gamma}a^{i+1}\otimes ... \otimes ^{\gamma}a^{n} = h_{i}(^{\gamma}a^{0},...,^{\gamma}a^{n})$,

showing that the homotopy h is compatible with the action of $\Gamma$ and it induces the homotopy of $inc^{\Gamma}_{*}tr^{\Gamma}_{*}$ to the identity.
This completes the proof of the theorem which extends well known results on Hochschild homology for $\Gamma$ acting trivially on A.

Besides $\Lambda = \mathds{Z}(G)$, the case $\Lambda = \mathds{Z}$ ($\Gamma$ acting trivially on $\mathds{Z}$) is also interesting. One means the consideration of the right $\Gamma$-equivariant derived functors $Ext^{n}_{\Lambda,\Gamma}(-,M)$ of the contravariant functor $Hom^{\Gamma}_{\Lambda}(-,M)$ from the category $\mathds{A}^{\Gamma}_{\Lambda}$ of $\Gamma$-equivariant left $\Lambda$-modules to the category of abelian groups, where $Hom^{\Gamma}_{\Lambda}(L,M)$ is the abelian group of $\Lambda$-homorphisms $f:L\rightarrow M$ compatible with the action of $\Gamma$. If $\Lambda = \mathds{Z}(G)$, we recover the functors $Ext^{n}_{\mathcal{P}}$ defined in [27] and $Ext^{n}_{{\mathds{Z}(G),\Gamma}}(\mathds{Z},L) \cong H^{n}_{\Gamma}(G,L), n\geq 0$.

\begin{rem}
Assume the group $\Gamma$ is acting on the $\kappa$-algebra A and consider the action of the group $\mathds{Z} \times \Gamma$ on the chain complex $C_{*}(A,A)$ as the composite of the action of $\mathds{Z}$ and the above defined action of $\Gamma$. Then the homology of the chain complex $(C_{*}(A.A))_{\mathds{Z} \times \Gamma}$ can be considered as the $\Gamma$-equivariant cyclic homology of the $\kappa$-algebra A for $\mathds{Q}\subset \kappa$.

\end{rem}

\;

\section{Extensions of crossed $\Gamma$-modules}

In this section the investigation of extensions of $\Gamma$-groups is continued for the class of $\Gamma$-groups endowed with a crossed $\Gamma$-module structure. These extensions are called $\Gamma$-extensions of crossed $\Gamma$-modules. As noted in the Introduction the extension theory of crossed modules  has been treated by many mathematicians. Our approach to extension theory of crossed modules substantially extends the class of relative extensions of epimorphisms of groups introduced and investigated by Loday [33].

A crossed $\Gamma$-module $(G,\mu)$ is a pair consisting of a $\Gamma$-group G and a $\Gamma$-homomorphism $\mu: G\rightarrow \Gamma$ ($\Gamma$ acting on itsel by conjugation) satisfying the Peiffer identity:

$^{\mu(g)}g' = gg'g_{-1}, g,g'\in G.$

A homomorphism from a crossed $\Gamma$-module $(G,\mu)$ to a crossed $\Gamma$-module $(G',\mu')$ is a $\Gamma$-homorphism $f:G\rightarrow G'$ such that $\mu'f = \mu$. Denote by $\mathrm{Cr}\Gamma$ the category of crossed $\Gamma$-modules. A crossed $\Gamma$-module $(G,\mu)$ will be called trivial if $\mu$ is the trivial map, $\mu(g) = e, g\in G$. There is an obvious equivalence between the category of trivial crossed $\Gamma$-modules and the category of $\Gamma$-modules.

A crossed $\Gamma$-module $(G,\mu)$ will be called elementary crossed $\Gamma$-module if $\mu$ is injective. It is equivalent to the inclusion crossed $\Gamma$-module $\overline{G}\overset{\overline{\mu}}\rightarrow \Gamma$, where $\overline{G} = \mu(G)$ is a normal subgroup of $\Gamma.$ If $(G,\mu)$ is a crossed $\Gamma$-module and A is a $\Gamma$-module satisfying the following property: $^{\gamma}a = a$, for $a\in A, \gamma\in Im\mu$, then A will be called crossed equivariant $\Gamma$-module. In particular $Ker \mu$ is a crossed equivariant $\Gamma$-module. It is obvious that any crossed equivariant $\Gamma$-module is a $\Gamma/Im \mu$-module.

\begin{defn}
Let
$(1)$ $0\rightarrow (A,1)\overset{\sigma}\rightarrow (X,\eta)\overset{\tau}\rightarrow (G,\mu)\rightarrow e$

be a sequence of crossed $\Gamma$-modules such that the induced sequence

$0\rightarrow A\overset{\sigma}\rightarrow X\overset{\tau}\rightarrow G\rightarrow e$

is an exact sequence of $\Gamma$-groups. Then the sequence (1) will be called $\Gamma$-extension of the crossed $\Gamma$-module $(G,\mu)$ by the crossed equivariant $\Gamma$-module A. In that case $\eta (X)$ acts trivially on A and $\sigma (A)$ belongs to the center of X.

If in addition there is a $\Gamma$-map $\alpha: (G,\mu)\rightarrow (X,\eta)$ such that the composite $\tau\alpha$ is the identity map, then it will be called $\Gamma$-extension with $\Gamma$-section map or $\Gamma$-equivariant extension of the crossed $\Gamma$-module $(G,\mu)$.
\end{defn}

Two $\Gamma$-extensions of $(G,\mu)$ by the crossed $\Gamma$-module (A,1)

$0\rightarrow (A,1)\overset{\sigma_{1}}\rightarrow (X_{1},\eta_{1})\overset{\tau_{1}}\rightarrow (G,\mu)\rightarrow e$,

and

$0\rightarrow (A,1)\overset{\sigma_{2}}\rightarrow (X_{2},\eta_{2})\overset{\tau_{2}}\rightarrow (G,\mu)\rightarrow e$

are called isomorphic if there is a $\Gamma$-homomorphism $\vartheta: (X_{1},\eta_{1})\rightarrow (X_{2},\eta_{2})$ inducing the identity map on (A,1) and $\tau_{2}\vartheta = \tau_{1}$.

Denote by $E((G,\mu),(A,1))$ and $E_{\Gamma}((G,\mu),(A,1))$ the set of isomorphism classes of $\Gamma$-extensions and of $\Gamma$-extensions having $\Gamma$-section map respectively. The contravariant functors $E( - ,(A,1))$, $E_{\Gamma}( - ,(A,1))$ on the category of crossed $\Gamma$-modules and the covariant functors  $E((G,\mu), - )$, $E_{\Gamma}((G,\mu), - )$ on the category of crossed equivariant $\Gamma$-modules to the category of sets are determined in a standard way. In particular, for the case of $E_{\Gamma}((G,\mu),(A,1))$ they are defined as follows:

Let $[E = (A,1)\overset{\sigma}\rightarrow (X,\eta)\overset{\tau}\rightarrow (G,\mu)]\in E_{\Gamma}((G,\mu),(A,1))$ with $\Gamma$-section map $\alpha$ and $f:(G',\mu')\rightarrow (G,\mu)$ be a $\Gamma$-homomorphism. By taking the fiber product $D = \{(x,g')\}$, $\tau(x) = f(g'), x\in X, g'\in G'$ of the diagram $X\overset{\sigma}\rightarrow G\overset{f}\leftarrow G'$ one obtains the $\Gamma$-extension $E' = (A,1)\overset{\sigma'}\rightarrow (D,\delta)\overset{p}\rightarrow (G',\mu')$, where $\sigma'(a) = (\sigma(a),e), p(x,g') = g', \delta(x,g') = \mu(x)$. The $\Gamma$-section map $\alpha':(G',\mu')\rightarrow (D,\delta)$ is given by $\alpha'(g') = (\alpha f(g'),g')$. This defines the contravariant functor $E((A,1), f): E_{\Gamma}((G,\mu),(A,1))\rightarrow E_{\Gamma}((G',\mu'),(A,1))$, $E((A,1),f)([E]) = [E']$. To define the covariant functor $E((G,\mu),(h)): E_{\Gamma}((G,\mu),(A,1))\rightarrow E_{\Gamma}((G,\mu),(A',1))$,  where $h:(A,1)\rightarrow (A',1)$ is a $\Gamma$-homomorphism, take the direct product $(A'\times X,\eta')$, $\eta' (a',x) = \eta (x)$, and the Cokernel ($ \beta$, $\eta''$) of the injection $\beta: A\rightarrow (A'\otimes X,\eta')$, $\beta(a) = (-h(a),\sigma (a))$, $\eta''([(a',x)]) = \eta'(x)$. This defines a $\Gamma$-extension $E'' = (A',1)\overset{\sigma''}\rightarrow (Cokernel \beta,\eta'')\overset{\tau''}\rightarrow (G,\mu)$, where $\sigma''(a') = [(a',1)]$, $\tau''[(a',x)] = \eta (x)$, with $\Gamma$-section map $\alpha'': (G,\mu)\rightarrow (Cokernel \beta, \eta''), \alpha'' = h'\alpha$, where $h'(x) = [(0,x)]$., and therefore the covariant functor $E((G,\mu),(h))([E]) = E''$.

To define the (co)homology and $\Gamma$-equivariant (co)homology of crossed $\Gamma$-modules two important classes will be defined in the category $\mathrm{Cr}\Gamma$ of crossed $\Gamma$-modules.

The objects of the first class $\mathfrak{P}_{\Gamma}$ of crossed $\Gamma$-modules are constructed as follows. Let $(G,\mu)$ be an arbitrary crossed $\Gamma$-module and take the free group $F(\Gamma\times G)$ generated by the couples $(\gamma,g)$, $\gamma \in \Gamma$, $g\in G$. There is an action of $\Gamma$ on $F(\Gamma\times G)$ given by $^{\gamma'}(\gamma,g) = (\gamma'\gamma,g)$, $\gamma, \gamma'\in \Gamma, g\in G$, and a $\Gamma$-homomorphism $\eta: F(\Gamma\times G)\rightarrow G, \eta (\gamma,g) = ^{\gamma}g$, inducing a $\Gamma$-homomorphism $\mu\eta: F(\Gamma\times G)\rightarrow \Gamma$. Consider the normal subgroup of  $F(\Gamma\times G)$ generated by the elements $^{\mu\eta}(x)x'x^{-1}x'^{-1}$ for $x,x'\in F(\Gamma\times G)$. Let $F_{(G,\mu)}$ denotes the quotient $\Gamma$-group $F(\Gamma\times G)/\{^{\mu\eta}(x)x'x^{-1}x'^{-1}\}$. Since $(G,\mu)$ is a crossed $\Gamma$-module, the $\Gamma$-homomorphism $\eta$ sends the normal subgroup $\{^{\mu\eta}(x)x'x^{-1}x'^{-1}\}$ to the unit. This yields a $\Gamma$-homomorphism $\eta':F_{(G,\mu)}\rightarrow G$ and a crossed $\Gamma$-module $(F_{(G,\mu)}, \mu\eta')$ which is called free crossed $\Gamma$-module generated by $(G,\mu)$ implying the canonical surjective homomorphism $\eta':(F_{(G,\mu)}, \mu\eta')\rightarrow (G,\mu)$. This construction was used by Loday to show the existence of the universal central relative extension of a group epimorphism [33]. The objects of the class $\mathfrak{P}_{\Gamma}$ are retracts of free crossed $\Gamma$-modules and are called projective crossed $\Gamma$-modules.

The construction of the second class $\mathfrak{P}_{\Gamma-e}$ of crossed $\Gamma$-modules is realized similarly. Consider the free group F(G) generated by the elements g , $g\in G$. There is an action of $\Gamma$ on F(G) given by $^{\gamma}|g| = |^{\gamma}g|$, $g\in G, \gamma\in \Gamma$, and let $\varphi: F(G)\rightarrow G$ be the canonical $\Gamma$-homomorphism, $\varphi(|g|) = g$, having a $\Gamma$-section map $\sigma:G\rightarrow F(G), \sigma(g) = |g|, g\in G$. This yields a crossed $\Gamma$-module $(F(G),\mu\varphi)$ and a homomorphism $\varphi: (F(G),\mu\varphi)\rightarrow (G,\mu)$ having a $\Gamma$-section map. The quotient $F_{\Gamma(G,\mu)} = F(G)/\{^{\mu\varphi}(x)x'x^{-1}x'^{-1}\}, x,x'\in F(G)$ provides a crossed $\Gamma$-module $(F_{\Gamma(G,\mu)}, \mu\varphi')$ induced by $\varphi$ which will be called $\Gamma$-equivariant free crossed $\Gamma$-module. and the canonical surjection $(F_{\Gamma(G,\mu)}, \mu\varphi')\rightarrow (G,\mu)$ having a $\Gamma$-section map. The class $\mathfrak{F}_{\Gamma-e}$ is consisting of all $\Gamma$-equivarianr free crossed $\Gamma$-modules.  The objects of the class $\mathfrak{P}_{\Gamma-e}$ are retracts of free crossed $\Gamma$-modules and are called $\Gamma$-equivariant projective crossed $\Gamma$-modules.

\begin{prop}
The classes $\mathfrak{P}_{\Gamma}$ and $\mathfrak{P}_{\Gamma-e}$ are projective classes in the category $\mathrm{Cr}\Gamma$ of crossed $\Gamma$-modules.
\end{prop}

$Proof$. To prove the class $\mathfrak{P}_{\Gamma}$ is projective it suffices to show that for any surjective homomorphism $f: (G',\mu')\rightarrow (G,\mu)$ of crossed $\Gamma$-modules and any homomorphism $h:F_{(L,\nu)}\rightarrow (G,\mu)$, where $F_{(L,\nu)}$ is a free crossed $\Gamma$-module, there is a homomorphism $h':F_{(L,\nu)}\rightarrow (G',\mu')$ that fh' = h.
For every element (e,l) of $\Gamma \times L, l\in L$, choose an element g' of G' such that f(g') = h([(e,l)]) and define the $\Gamma$-map $\Gamma \times L\rightarrow G'$ sending $(\gamma, l)$ to $^{\gamma}g'$ which induces the required homomorphism h'.

For the class $\mathfrak{P}_{\Gamma-e}$ it suffices to show that for any surjective homomorphism $f: (G',\mu')\rightarrow (G,\mu)$ of crossed $\Gamma$-modules having a $\Gamma$-section map and any homomorphism $h:F_{\Gamma(L,\nu)}\rightarrow (G,\mu)$, where $F_{\Gamma(L,\nu)}$ is a $\Gamma$-equivariant free crossed $\Gamma$-module, there is a homomorphism $h':F_{\Gamma(L,\nu)}\rightarrow (G',\mu')$ that fh' = h. It is easily checked that h' can be defined as $h'([x]) = \sigma h([x])$ where $x\in F(L)$ an $\sigma$ is the section-map of f.
This completes the proof of the Proposition.

There is a $\Gamma$-homomorphism $F(\Gamma \times G)\rightarrow F(G)$ sending the generator $|(\gamma,g)|$ to $|^{\gamma}g|$ inducing $\Gamma$-homomorphism $F_{(G,\mu)}\rightarrow F_{\Gamma(G,\mu)}$ and therefore a natural morphism $\omega:\mathfrak{P}_{\Gamma}\rightarrow
\mathfrak{P}_{\Gamma-e}$ from the projective class  $\mathfrak{P}_{\Gamma}$ to the projective class $\mathfrak{P}_{\Gamma-e}$.

For the cohomological interpretation of the abelian group of $\Gamma$-extensions of crossed $\Gamma$-modules the right derived functors $R_{\mathfrak{\bar{P}}}^{n}T, n\geq 0$, of a contravariant functor T from the category $\mathfrak{A}$ with finite inverse limits to the category of abelian groups with respect to a projective class $\mathfrak{\bar{P}}$ will be defined. The case of left derived functors of a covariant functor to the category of abelian groups or to the category of groups was considered in [47] and [24,26] respectively.

To this aim let us recall some definitions given in [26].

\begin{defn}

 A $\mathfrak{\bar{P}}$-projective resolution of an object A of the category $\mathfrak{A}$ is a pseudo-simplicial projective object over A, $P_{*}\rightarrow A$, which is $\mathfrak{\bar{P}}$-exact and $\mathfrak{\bar{P}}$-epimorphic.

\end{defn}

Since the category $\mathfrak{A}$ contains finite inverse limits, every object A admits a $\mathfrak{\bar{P}}$-projective resolution which is unique up to simplicial homotopy.

\begin{defn}
The right derived functors $R_{\mathfrak{\bar{P}}}^{n}T$ of the contravariant functor T with respect to the projective class $\mathfrak{\bar{P}}$ are given by $R_{\mathfrak{\bar{P}}}^{n}T(A) = H_{n}T(P_{*}), n\geq 0.$
\end{defn}
 It is obvious that the category $\mathrm{Cr}\Gamma$ of crossed $\Gamma$-modules is a category with finite inverse limits. Denote by $Hom_{\Gamma}((G,\mu),(A,1))$ the abelian group of $\Gamma$-homomorphisms from the $\Gamma$-group G to the $\Gamma$-module A and consider the contravariant functor $Hom_{\Gamma}(-,(A,1))$ from the category $\mathrm{Cr}\Gamma$ to the category of abelian groups.

 The cohomology and the $\Gamma$-equivariant cohomology of crossed $\Gamma$-modules will now be defined by using the right derived functors of the functor $Hom_{\Gamma}(-,(A,1))$ with respect to the projective classes $\mathfrak{P}$ and $\mathfrak{P}_{\Gamma}$ respectively. Namely

\begin{defn}

The n-th cohomology and $\Gamma$-equivariant cohomology of the crossed $\Gamma$-module $(G,\mu)$ with coefficients in a $\Gamma$-module A are given by

$H^{n}_{\mathfrak{P}_{\Gamma}}((G,\mu),A) = R^{n-1}_{\mathfrak{P}_{\Gamma}}Hom_{\Gamma}((G,\mu),(A,1))$ and

$H^{n}_{\mathfrak{P}_{\Gamma-e}}((G,\mu),A) = R^{n-1}_{\mathfrak{P}_{\Gamma-e}}Hom_{\Gamma}((G,\mu),(A,1))$, $n\geq 1$,

respectively

\end{defn}

\begin{thm}

We have

$1) H^{1}_{\mathfrak{P}_{\Gamma}}((G,\mu),A) =  H^{1}_{\mathfrak{P}_{\Gamma-e}}((G,\mu),A) \cong Hom_{\Gamma}((G,\mu),(A,1)).$

$2) H^{2}_{\mathfrak{P}_{\Gamma}}((G,\mu),A) \cong E((G,\mu),(A,1))$  and

$H^{2}_{\mathfrak{P}_{\Gamma-e}}((G,\mu),A) \cong E_{\Gamma}((G,\mu),(A,1)),$

where A is crossed equivariant $\Gamma$-module.

$3) H^{2}_{\mathfrak{P}_{\Gamma}}((Im \mu,\sigma),A) \cong \mathcal{E}xt(\Gamma/Im \mu,\Gamma,A)$,

where $\mathcal{E}xt(\Gamma/Im \mu,\Gamma,A)$ is the abelian group of relative extensions of $(\Gamma/Im\mu,\Gamma)$ defined by Loday [33] and $(Im \mu, \sigma,)$ is the induced inclusion crossed $\Gamma$-module, $\sigma: Im \mu\hookrightarrow \Gamma$.

4) The short exact sequence of $\Gamma$-modules

\begin{equation}
 0\rightarrow A'\overset{f'}\rightarrow A\overset{f}\rightarrow A''\rightarrow 0
\end{equation}

gives rise a long exact cohomology sequence for the $\Gamma$-module $(G,\mu)$ if the $\Gamma$-modules of (6.1) are crossed equivariant:

$0\rightarrow Hom_{\Gamma}((G,\mu), A')\rightarrow Hom_{\Gamma}((G,\mu), A))\rightarrow Hom_{\Gamma}((G,\mu), A''))\rightarrow H^{2}_{\mathfrak{P}_{\Gamma}}((G,\mu),A')\rightarrow H^{2}_{\mathfrak{P}_{\Gamma}}((G,\mu),A)\rightarrow H^{2}_{\mathfrak{P}_{\Gamma}}((G,\mu),A'') \rightarrow H^{3}_{\mathfrak{P}_{\Gamma}}((G,\mu),A')\rightarrow ... \rightarrow H^{n}_{\mathfrak{P}_{\Gamma}}((G,\mu),A')\rightarrow H^{n}_{\mathfrak{P}_{\Gamma}}((G,\mu),A)\rightarrow H^{2}_{\mathfrak{P}_{\Gamma}}((G,\mu),A'')\rightarrow H^{n+1}_{\mathfrak{P}_{\Gamma}}((G,\mu),A')\rightarrow ...$.

If the sequence (6.1) of $\Gamma$-modules possesses the $\Gamma$-property then it induces the long cohomology sequence

$0\rightarrow Hom_{\Gamma}((G,\mu), A')\rightarrow Hom_{\Gamma}((G,\mu), A))\rightarrow Hom_{\Gamma}((G,\mu), A''))\rightarrow H^{2}_{\mathfrak{P}_{\Gamma-e}}((G,\mu),A')\rightarrow H^{2}_{\mathfrak{P}_{\Gamma-e}}((G,\mu),A)\rightarrow H^{2}_{\mathfrak{P}_{\Gamma-e}}((G,\mu),A'') \rightarrow H^{3}_{\mathfrak{P}_{\Gamma-e}}((G,\mu),A')\rightarrow ... \rightarrow H^{n}_{\mathfrak{P}_{\Gamma-e}}((G,\mu),A')\rightarrow H^{n}_{\mathfrak{P}_{\Gamma-e}}((G,\mu),A)\rightarrow H^{2}_{\mathfrak{P}_{\Gamma-e}}((G,\mu),A'')\rightarrow H^{n+1}_{\mathfrak{P}_{\Gamma-e}}((G,\mu),A')\rightarrow ...$.
\end{thm}

$Proof$. 1) Obvious.

2) Consider the canonical $\mathfrak{P}_{\Gamma}$-projective resolution of $(G,\mu)$

 $...\overset{\tau_{2}}\rightarrow (K_{2},l_{2})\underset{\lambda^{2}_{2}}{\overset{\lambda^{0}_{2}}\Rrightarrow} P_{((K_{1},l_{1}),\mu_{1})}\overset{\tau_{1}}\rightarrow (K_{1},l_{1})\underset{\lambda^{1}_{1}}{\overset{\lambda^{0}_{1}}\rightrightarrows} P_{((G,\mu),\mu \eta')}\overset{\tau_{0}}\rightarrow(G,\mu)$,

 where $(K_{1},l_{1})$ is the simplicial kernel of $\tau_{1}$ and $(K_{2},l_{2})$ is the simplicial kernel of $(\tau_{1}\lambda^{0}_{1},\tau_{1}\lambda^{1}_{1})$. Let $f: P_{(K_{1}n ,l_{1})}\rightarrow A$ be a $\Gamma$-homomorphism such that $\underset{i=0}{\overset{2}\sum}\tau_{2}\lambda^{i}_{2}f = o$. Therefore $\underset{i=0}{\overset{2}\sum}\lambda^{i}_{2}f = o$ and this implies a $\Gamma$-homomorphism $f_{1}: (K_{1},l_{1})\rightarrow A$ given by $f_{1}(x) = f(y)$ for $y = \tau_{1}(x), x\in (K_{1},l_{1})$. The correctness follows from the fact if $y_{1} = \tau_{1}(x)$ then the triple $(y_{1}y^{-1},y_{1}y^{-1},y_{1}y^{-1})$ belongs to $(K_{2},l_{2})$ implying $f(y_{1}y^{-1}) = 0$. By the  same argument we have $f(x,x) = 0$ for any $x\in P_{((G,\mu),\mu \eta')}$.

Now by using the diagram

$(K_{1},l_{1})\underset{\lambda^{1}_{1}}{\overset{\lambda^{0}_{1}}\rightrightarrows} P_{((G,\mu),\mu \eta')}\overset{\tau_{0}}\rightarrow(G,\mu)$,

and the $\Gamma$-homomorphism $f_{1}: (K_{1},l_{1})\rightarrow A$ the following crossed $\Gamma$-extension of $(G,\mu)$ is constructed. Take the direct product $A\times P_{(G,\mu)}$ with component wise action of $\Gamma$ on it. We obtain the crossed $\Gamma$-module $(A\times P_{(G,\mu)}, \bar{\mu})$ where $\bar{\mu}([a,x]) = \mu\eta'(x)$, $x\in P_{(G,\mu)}, a\in A$. By introducing on $A\times P_{(G,\mu)}$ the following equivalence relation:

$(a,x) \sim (b,y)$ if $\tau_{0}(x) = \tau_{0}(y)$ and $a\cdot f_{1}(x,y) = b$, this yields the crossed $\Gamma$-module $(A\times P_{(G,\mu)}/ \sim, \alpha) $, $\alpha([a,x]) = \mu\eta'(x)$, and the mentioned crossed $\Gamma$-extension of $(G,\mu)$:

$E = 0\rightarrow (A,1)\overset{\sigma}\rightarrow (A\times P_{(G,\mu)}/ \sim, \alpha)\overset{\beta}\rightarrow (G,\mu)\rightarrow e$, where $\sigma(a) = [ (a,e)], \beta([a,x]) = \tau_{0}(x)$. This allows to define correctly a homomorphism $\vartheta:H^{2}_{\mathfrak{P}}((G,\mu),A)\rightarrow E((G,\mu),(A,1))$ sending $[f]$ to $[E]$. Conversely for any crossed $\Gamma$-extension E of $(G,\mu)$:

$E = 0\rightarrow (A,1)\rightarrow (X,\eta)\rightarrow (G,\mu)\rightarrow e$ the $\Gamma$-homomorphism $\tau_{0}: (P_{(G,\mu)},\mu\eta')\rightarrow (G,\mu)$ induces a $\Gamma$-homomorphism $h': P_{(G,\mu)},\mu\eta')\rightarrow (X,\eta)$ and its composite with the homomorphism $(X,\eta)\rightarrow (G,\mu)$ is equal to $\tau_{0}$. Therefore this implies a homomorphism $g': K_{1},l_{1})\rightarrow (A,1)$ such that the composite of g' with the homomorphism $(A,1)\rightarrow (X,\eta)$ is equal to $\gamma^{0}_{1}(\gamma^{1}_{1})^{-1}$. The $\Gamma$-homomorphism $g = g'\tau_{1}:  P_{((G,\mu),\mu \eta')}\rightarrow A$ satisfies the equality $\underset{i=0}{\overset{2}\sum}\tau_{2}\lambda^{i}_{2}g = 0$ and implies the homomorphism  $\vartheta':E((G,\mu),(A,1))\rightarrow H^{2}_{\mathfrak{P}_{\Gamma}}((G,\mu),A)$ sending $[E]$ to $[g]$ such that the homomorphisms $\vartheta$ and $\vartheta'$ are inverse to each other.

3) First let us recall the definition of relative extensions of group epimorphisms.

\begin{defn}$[33]$

Let $\nu: N\rightarrow Q$ be an epimorphism of groups. A relative extension of (Q,N) is given by an exact sequence of groups

$1\rightarrow L\overset{\lambda}\rightarrow M\overset{\mu}\rightarrow N\overset{\nu}\rightarrow Q\rightarrow 1$

and an action $\eta$ of N on M such that $(M,\mu)$ is a crossed N-module.
\end{defn}

It is evident that every relative extension $(M,\mu)$ of (Q,N) induces the N-extension of $(Im\mu,\sigma)$

$$
0\rightarrow (L,1)\rightarrow (M,\mu)\rightarrow ((Im\mu,\sigma)\rightarrow 1
$$

where $\sigma: Im\mu\hookrightarrow N$ is the inclusion crossed N-module and one gets a map $\mathcal{E}xt(Q,N,L)\rightarrow E((Im\mu,\sigma),(L,1))$.
Conversely, any N-extension of $(Im\mu, \sigma)$

$0\rightarrow (L,1)\rightarrow (X,\mu)\rightarrow (Im\mu,\sigma)\rightarrow 1$

induces a relative extension of $(N/Im\mu,N)$

$1\rightarrow L\rightarrow X\rightarrow N\rightarrow N/Im\mu\rightarrow 1$

and therefore a map $E((Im\mu,\sigma),(L,1))\rightarrow \mathcal{E}xt(N/Im\mu,N,L) \cong \mathcal{E}xt(Q,N,L)$ which is the inverse of the map $\mathcal{E}xt(Q,N,L)\rightarrow E((Im\mu,\sigma),(L,1))$. It remains to apply 2).

4)Let

$ 0\rightarrow A'\overset{f'}\rightarrow A\overset{f}\rightarrow A''\rightarrow 0 $

be a short exact sequence of crossed equivariant $\Gamma$-modules. It suffices to prove the surjection of the homomorphism $Hom_{\Gamma}((P_{(G,\mu)}, \mu\eta'),(A,1)))\rightarrow Hom_{\Gamma}((P_{(G,\mu)}, \mu\eta'),(A'',1))$ induced by the $\Gamma$-homomorphism f. Let $h: P_{(G,\mu)}\rightarrow A''$ be a $\Gamma$-homomorphism and $\tau: P(\Gamma\times G)\rightarrow P_{(G,\mu)}$ be the canonical surjection. For $h\tau(e,g), g\in G$, take $a\in A$ such that $f(a)= h\tau(e,g)$. This yields the $\Gamma$-homomorphism $\bar{h}:P(\Gamma\times G)\rightarrow A$ given by $\bar{h}(\gamma,g) = ^{\gamma} a$, $\gamma\in \Gamma, g\in G$ such that $f\bar{h} = h\tau$. It is clear that the homomorphism $\bar{h}$ sends every element of the subgroup $\{^{\mu\eta}(x)x'x^{-1}x'^{-1}\}$  to the unit, since A is a crossed equivariant $\Gamma$-module and we obtain a $\Gamma$-homomorphism $h':P_{(G,\mu)}\rightarrow A$ induced by $\bar{h}$ such that $fh' = h$.

If the sequence (6.1) possesses the $\Gamma$-property, then the $\Gamma$-homomorphism $h:P_{(G,\mu)}\rightarrow A''$ induces  the $\Gamma$-homomorphism $h'': P_{\Gamma(G,\mu)}\rightarrow A''$ which implies the $\Gamma$-homomorphism $h': P_{\Gamma(G,\mu)}\rightarrow A$ given by $h'(x) = \sigma h''(x)$, $x \in G$, where $\sigma$ is the $\Gamma$-section map of (6.1). This yields the equality $fh'({^{\mu\eta}(x)x'x^{-1}x'^{-1}}) = f(^{\mu\eta}(h'(x)(h'(x)^{-1}) = 1$. The $\Gamma$-property of (6.1) implies the equality $h'(^{\mu\eta}(x)x^{-1}) = 1$. Therefore the subgroup $\{^{\mu\eta}(x)x'x^{-1}x'^{-1}\}$ goes to unit by the $\Gamma$-homomorphism h'. It is evident that the $\Gamma$-homomorphism $\bar{h}: P_{\Gamma(G,\mu)}\rightarrow A$ induced by h' satisfies the equality $f\bar{h} = h.$ This completes the proof of the theorem.

\;

\section{Homology and central $\Gamma$-extensions of crossed $\Gamma$-modules}

To define the homology of crossed $\Gamma$-modules the left derived functors $L^{\mathfrak{P}}_{n}T$, $n\in 0$, of a covariant functor $T:\mathfrak{A}\rightarrow Gr$ from the category $\mathfrak{A}$ with finite inverse limits to the category $\mathrm{Gr}$ of groups will be used [24,26].

\begin{defn} $[26]$

Let $P_{*}\rightarrow A$ be a pseudo-simplicial projective resolution over the object A of the category $\mathfrak{A}$:

 $P_{*} = ...\overset{\tau_{2}}\rightarrow P_{2}\underset{\lambda^{2}_{2}}{\overset{\lambda^{0}_{2}}\Rrightarrow}  P_{1}\underset{\lambda^{1}_{1}}{\overset{\lambda^{0}_{1}}\rightrightarrows}  P_{0}\overset{\tau_{0}}\rightarrow A$,

 and consider the chain complex $(L_{*}T(P_{*}), d_{*})$, where

 $L_{n}T(P_{*}) = T(P_{n})\cap KerT(\lambda^{0}_{n})\cap ... \cap KerT(\lambda^{n}_{n-1})$ , $n\geq 0$, and $d_{n}: L_{n}T(P_{*})\rightarrow L_{n-1}T(P_{*})$  is the restriction of $T(\lambda^{n}_{n})$ on $L_{n}T(P_{*})$.

 The n-th homology group of the chain complex $(L_{*}T(P_{*}), d_{*})$ defines the n-th left derived functor $L^{\mathfrak{P}_{\Gamma}}_{n}T$ of T with respect to the projective class $\mathfrak{P}_{\Gamma}$.
\end{defn}

If the values of the functor T belong to the category of abelian groups one can also use the definition of Tierney - Vogel [47] by considering the homology groups of the chain complex $JT(P_{*})$:

$JT(P_{*}) = \{T(P_{n}), \delta_{n}, n\geq 0\}$, where $\delta_{n} = \underset{i=0}{\overset{n}\sum}(-1)^{i}T(\delta^{n}_{i})$.

The natural homomorphism $LT(P_{*})\rightarrow JT(P_{*})$ induces an isomorphism of their homology groups and the proof of this assertion is completely similar to the proof for simplicial groups [35].

\begin{defn}

The n-th homology group of the crossed $\Gamma$-module $(G,\mu)$ with coefficients in the $\Gamma$-module A and with respect to the projective class $\mathfrak{P}_{\Gamma}$ is given by

$H^{\mathfrak{P}_{\Gamma}}_{n}((G,\mu),A) = L^{\mathfrak{P}_{\Gamma}}_{n-1}(I(G)\otimes_{G\rtimes \Gamma}A)$, $n\geq 1$.

The $\Gamma$-equivariant n-th homology group of the crossed $\Gamma$-module $(G,\mu)$ with coefficients in the $\Gamma$-module A and with respect to the projective class $\mathfrak{P}_{\Gamma-e}$ is defined as

$H^{\mathfrak{P}_{\Gamma-e}}_{n}((G,\mu),A) = L^{\mathfrak{P}_{\Gamma-e}}_{n-1}(I(G)\otimes_{G\rtimes \Gamma}A)$, $n\geq 1$.

\end{defn}

\begin{prop}

One has

1) $H^{\mathfrak{P}_{\Gamma}}_{n}((G,\mu),A) = 0$ for $n\geq 2$ and $H^{\mathfrak{P}_{\Gamma}}_{1}((G,\mu),A) = I(G)\otimes_{G\propto \Gamma}A$, if $(G,\mu)$ is a projective crossed $\Gamma$-module.

2)If $\Gamma$ acts trivially on A, then $H^{\mathfrak{P}_{\Gamma}}_{1}((G,\mu),A) = G/[G,G]_{\Gamma} \otimes A.$

\end{prop}

$Proof$. 1) Let $(P,\mu)$ be a projective crossed $\Gamma$-module and $((Y_{*},\delta_{*}),\tau, (P,\mu))$ be a $\mathfrak{P}_{\Gamma}$-projective resolution of $(P,\mu)$. Since P is projective, there is a $\Gamma$-homomorphism $h:P\rightarrow Y_{0}$ such that $\tau h = 1$ and inducing the left contractibility $h_{n}:Y_{n}\rightarrow Y_{n+1}$ , $n\geq 0$. It follows that the abelian augmented pseudo-simplicial group $(I(Y_{*})\otimes_{G\rtimes \Gamma}A), I(\tau)\otimes_{G\rtimes \Gamma}A), I(P)\otimes_{G\rtimes \Gamma}A)$ is left contractive and therefore aspherical. This also yields $L^{\mathfrak{P}_{\Gamma}}_{0}(I(G)\otimes_{G\rtimes \Gamma}A) = I(G)\otimes_{G\rtimes \Gamma} A$.

2) First it will be shown that every exact sequence of $\Gamma$-groups

$G_{1}\overset{\tau_{1}}\rightarrow G_{2}\overset{\tau_{2}}\rightarrow G\rightarrow e$.

 induces the exact sequence

$H_{1}^{\Gamma}(G_{1})\overset{H_{1}^{\Gamma}(\tau_{1})}\rightarrow H_{1}^{\Gamma}(G_{2})\overset{H_{1}^{\Gamma}(\tau_{2})}\rightarrow H_{1}^{\Gamma}(G)\rightarrow 0$.

By using the exact sequence [27]

$0\rightarrow \Gamma G/[G,G]\cap\Gamma G\rightarrow H_{1}(G)\rightarrow H^{\Gamma}_{1}(G)\rightarrow 0$

one gets the following commutative diagram with exact rows and columns

$$
\xymatrix{
 0 \ar[d] & 0 \ar[d] & 0 \ar[d] \\
Ker \xi \ar[d]\ar[r] &   \Gamma G_{2}/[G_{2},G_{2}]\cap\Gamma G_{2} \ar[d]\ar[r]^{\xi} & \ar[d]\ar[r]   \Gamma G/[G,G]\cap\Gamma G \ar[d]\ar[r] & 0\\
KerH_{1}(\tau_{2}) \ar[d]^{\sigma}\ar[r] & H_{1}(G_{2}) \ar[d]\ar[r]^{H_{1}(\tau_{2})}  & H_{1}(G) \ar[d]\ar[r] & 0\\
KerH^{\Gamma}_{1}(\tau_{2})\ar[r] & H^{\Gamma}_{1}(G_{2}) \ar[d]\ar[r]^{H^{\Gamma}_{1}(\tau_{2})} & H^{\Gamma}_{1}(G)\ar[d]\ar[r] & 0\\
 & 0 & 0 &
}
$$

showing that the homomorphism $\sigma$ is surjective. Consider now the commutative diagram

$$
\xymatrix{
 H_{1}(G_{1}) \ar[d]^{H_{1}(\tau_{1})}\ar[r]^{\delta} & KerH(\tau_{2})\ar[d]^{\sigma}\\
H^{\Gamma}_{1}(G_{1}) \ar [r]^{\delta'} & KerH^{\Gamma}_{1}(\tau_{2}).\\
 }
$$

By Stallings - Stammbach exact homology sequence [44] one concludes that $\delta$ is surjective implying the surjection of $\delta'$ and therefore the exactness of the required sequence. The isomorphism $I(G)\otimes_{G\propto \Gamma}A \cong H^{\Gamma}_{1}(G)\otimes A$ [27] gives us the exactness of the sequence

$I(G_{1})\otimes_{G_{1}\rtimes \Gamma}A\rightarrow I(G_{2})\otimes_{G_{2}\rtimes \Gamma}A\rightarrow I(G)\otimes_{G\rtimes \Gamma}A\rightarrow 0$.

It is now easily checked that if $((P_{*}, \mu_{*}), \tau , (G,\mu))$ is a $\mathfrak{P}_{\Gamma}$-projective resolution of the crossed $\Gamma$-module  $(G,\mu)$ then $\pi_{0}(I(P_{*})\otimes_{P_{*}\rtimes \Gamma}A)$ is isomorphic to $I(G)\otimes_{G\rtimes \Gamma}A \cong G/[G,G]_{\Gamma}\otimes A$. This completes the proof of the proposition.

\begin{rem}

The results of Prop.7.3 also hold for the $\Gamma$-equivariant homology $H^{\mathfrak{P}_{\Gamma-e}}_{*}((G,\mu),A)$.

\end{rem}

In what follows $H^{\mathfrak{P}_{\Gamma}}_{n}((G,\mu),Z)$ and $H^{\mathfrak{P}_{\Gamma-e}}_{n}((G,\mu),Z)$  are denoted $H^{\mathfrak{P}_{\Gamma}}_{n}(G,\mu)$ and $H^{\mathfrak{P}_{\Gamma-e}}_{n}(G,\mu)$ respectively.

\begin{defn}

The $\Gamma$-extension of the crossed $\Gamma$-module $(G,\mu)$ by a $\Gamma$-module A

 $0\rightarrow (A,1)\overset{\sigma}\rightarrow (U,\eta)\overset{\tau}\rightarrow (G,\mu)\rightarrow e$

is called central if $\Gamma$ acts trivially on A. It is called universal if for any central $\Gamma$-extension $(Y,\delta)$ of $(G,\mu)$

 $0\rightarrow (B,1)\rightarrow (Y,\delta)\rightarrow (G,\mu)\rightarrow e$

there is a unique $\Gamma$- homomorphism $(U,\eta)\rightarrow (Y,\delta)$ over $(G,\mu)$.

The $\Gamma$-extension of the crossed $\Gamma$-module $(G,\mu)$ by a $\Gamma$-module A with $\Gamma$-section map

 $0\rightarrow (A,1)\overset{\sigma}\rightarrow (U,\eta)\overset{\tau}\rightarrow (G,\mu)\rightarrow e$

is called $\Gamma$-equivariant central if $\Gamma$ acts trivially on A. It is called $\Gamma$-equivariant universal if for any central crossed $\Gamma$-extension  $(Y,\delta)$ of $(G,\mu)$ with $\Gamma$-section map

 $0\rightarrow (B,1)\rightarrow (Y,\delta)\rightarrow (G,\mu)\rightarrow e$

there is a unique $\Gamma$- homomorphism $(U,\eta)\rightarrow (Y,\delta)$ over $(G,\mu)$.

For the construction of the universal $\Gamma$-extension and the $\Gamma$-equivariant universal $\Gamma$-extension of the crossed $\Gamma$-module $(G,\mu)$ the projective classes  $\mathfrak{P}_{\Gamma}$ and $\mathfrak{P}_{\Gamma-e}$ will be used respectively.
\end{defn}

Consider the above mentioned free group $F(\Gamma,G)$ with $\Gamma$-action $^{\gamma'}(\gamma,g) = (\gamma'\gamma,g)$ and the $\Gamma$-homomorphism $\eta': F_{(G,\mu)}\rightarrow G$ given by $\eta'([(\gamma,g)]) = ^{\gamma}g$. Denote by R the kernel of $\eta'$ and by P the quotient of the free crossed $\Gamma$-module $F_{(G,\mu)}$ by the normal subgroup generated by the elements $[^{\gamma}r\cdot r^{-1}], r\in R, \gamma\in \Gamma$. This yields a crossed $\Gamma$-module $(P,\tau)$, where $\tau$ is induced by the canonical homomorphism $\eta':(F_{(G,\mu)}, \mu\eta')\rightarrow (G,\mu)$.

For the $\Gamma$-equivarianr case take the free group F(G) with $\Gamma$-action $^{\gamma}|g| = |^{\gamma}g| , g\in G.$ and the kernel L of the  homomorphism $\varphi':F_{\Gamma(G,\mu)}\rightarrow G$ induced by the canonical homomorphism $\varphi:F(G)\rightarrow G$.  Let $P_{\Gamma}$ be the quotient of the $\Gamma$-equivariant free crossed $\Gamma$-module $F_{\Gamma(G,\mu)}$ by the normal subgroup generated by the elements $[^{\gamma}l\cdot l^{-1}]$, $l\in L, \gamma\in \Gamma$. This yields a crossed $\Gamma$-module $(P_{\Gamma},\tau_{\Gamma})$, where $\tau_{\Gamma}$ is induced by the canonical surjection $F_{\Gamma(G,\mu)},\mu\varphi')\rightarrow (G,\mu)$.

\begin{defn}

The crossed $\Gamma$-module $(G,\mu)$ is called $\Gamma$-perfect if $H^{\mathfrak{P}_{\Gamma}}_{1}(G,\mu) = H^{\mathfrak{P}_{\Gamma-e}}_{1}(G,\mu) = 0.$

\end{defn}

\begin{thm}

1)  A central $\Gamma$-extension $(U,\eta)$ of a crossed $\Gamma$-module $(G,\mu)$ is universal if and only if it is $\Gamma$-perfect and every central  $\Gamma$-extension of $(U,\eta)$ splits.

2) If $(G,\mu)$ is $\Gamma$-perfect then the $\Gamma$-extension

 $0\rightarrow (R',1)\rightarrow ([P,P]_{\Gamma}/[P,R]_{\Gamma},\bar{\eta})\overset{\tau'}\rightarrow (G,\mu)\rightarrow e$

 is the universal $\Gamma$-extension of $(G,\mu)$, where $\tau'$ is induced by $\tau$ and $\bar{\eta}$ is induced by $\mu\eta'$.

 3) A central $\Gamma$-equivariant extension $(U,\eta)$ of a crossed $\Gamma$-module $(G,\mu)$ is $\Gamma$-equivariant universal if and
 only if it is $\Gamma$-perfect and every central $\Gamma$-equivariant extension of $(U,\eta)$ splits.

 4) If $(G,\mu)$ is $\Gamma$-perfect then the $\Gamma$-extension

$0\rightarrow (L',1)\rightarrow ([P_{\Gamma},P_{\Gamma}]_{\Gamma}/[P_{\Gamma},L]_{\Gamma}, \bar{\varphi})\overset{\tau'}\rightarrow (G,\mu)\rightarrow e$

 is the $\Gamma$-equivariant universal $\Gamma$-extension of $(G,\mu)$, where $\tau'$ is induced by $\tau$ and $\bar{\varphi}$ is induced by $\mu\varphi$.

\end{thm}

$Proof$. 1)The way follows to the classical proof for central group extensions [37] and that has been also realized for central $\Gamma$-equivariant group extensions [27]. Let

$E = 0\rightarrow (C,1)\overset{\alpha}\rightarrow (U,\eta)\overset{\beta}\rightarrow (G,\mu)\rightarrow e$

be the universal $\Gamma$-extension of $(G,\mu)$. If $(U,\eta)$ is not perfect, there is two distinct morphisms $f_{1}, f_{2}: (U,\eta)\rightarrow (U/[U,U]_{\Gamma}, \eta')$ from E to $F = 0\rightarrow (U/[U,U]_{\Gamma},1)\overset{\sigma}\rightarrow (U/[U,U]_{\Gamma}, \eta')\overset{\tau}\rightarrow (G,\mu)\rightarrow e$ over $(G,\mu)$, where $f_{1}(x) = (1, \beta(x))$, $f_{2}(x) = (\psi(x),\beta(x))$, $x\in U$, and $\psi: U\rightarrow U/[U,U]_{\Gamma}$ is the canonical homomorphism. That is in contradiction with the universality of E.

Let

$F = 0\rightarrow (D,1)\overset{\sigma}\rightarrow (W,\delta)\overset{\omega}\rightarrow (U,\eta)\rightarrow e$

be a central $\Gamma$-extension of $(U,\eta)$. The sequence of $\Gamma$-modules

$0\rightarrow (Ker \beta\omega,1)\rightarrow (W,\delta)\rightarrow (G,\mu)\rightarrow e$ is a central $\Gamma$-extension of $(G,\mu)$. The universality of E yields a homomorphism $f':(U,\eta)\rightarrow (W,\beta)$ over $(G,\mu)$ and the composite $\omega f'$ is also a homomorphism over $(G,\mu)$ implying the equality $\omega f' = 1$ and the splitting of F.

2)Let

$0\rightarrow (D,1)\rightarrow (X,\delta)\rightarrow \overset{\omega}\rightarrow (G,\mu)\rightarrow 1$

be a central $\Gamma$-extension of $(G,\mu)$. Then there is a homomorphism $f: (P,\tau)\rightarrow (X,\delta)$ of crossed $\Gamma$-modules over $(G,\mu)$ such that the diagram

$$
\xymatrix{
([P,P]_{\Gamma}/[P,R]_{\Gamma}, \eta')\ar[r]^(0.7){\tau'}  \ar[d]^{f'}   &(G,\mu) \ar[d]^\parallel\\
(X,\delta) \ar [r]^{\omega} &(G,\mu).\\,
 }
$$

where $\tau'$, $\eta'$ and $f'$ are induced by $\tau$, $\mu$ and $f$ respectively. The homomorphism $\tau'$ is surjective, since $(G,\mu)$ is $\Gamma$-perfect and therefore $\tau'$ is a central $\Gamma$-extension of $(G,\mu)$. Every $\Gamma$- homomorphism $[P,P]_{\Gamma}/[P,R]_{\Gamma} \rightarrow Ker \omega$ is trivial implying the uniqueness of $f'$ over $(G,\mu)$.

3)and 4) We omit the proof since it goes along the same lines as the proof of 1) and 2) respectively by replacing in particular the crossed $\Gamma$-module $(P,\tau)$ by the crossed $\Gamma$-module $(P_{\Gamma},\tau_{\Gamma})$.

$Conclusion$ A crossed $\Gamma$-module has the universal $\Gamma$-extension and the $\Gamma$-equivariant universal $\Gamma$-extension if and only if it is $\Gamma$-perfect.

\begin{thm}
1) Let

$0\rightarrow (A,1)\rightarrow (B,\delta)\overset{\vartheta}\rightarrow (G,\mu)\rightarrow 1$

be a $\Gamma$-extension of $(G,\mu)$ and $\tau: (F,\eta)\rightarrow (B,\delta)$ be a free presentation of $(B,\delta)$. Then there is an exact sequence

$0\rightarrow U\rightarrow H^{\mathfrak{P}_{\Gamma}}_{2}(B,\delta)\rightarrow H^{\mathfrak{P}_{\Gamma}}_{2}(G,\mu)\overset{\rho}\rightarrow A/[B,A]_{\Gamma}\rightarrow H^{\mathfrak{P}_{\Gamma}}_{1}(B,\delta)\rightarrow H^{\mathfrak{P}_{\Gamma}}_{1}(G,\mu)\rightarrow 0$,

where U is the kernel of  $[F,S]_{\Gamma}/[F,R]_{\Gamma}\rightarrow [B,A]_{\Gamma}$, $R = Ker \tau $ and $S = Ker \theta \tau$.

2) If

$0\rightarrow (A,1)\rightarrow (B,\delta)\overset{\vartheta}\rightarrow (G,\mu)\rightarrow 1$

is a $\Gamma$-equivariant extension of $(G,\mu)$ and $\tau_{\Gamma}: (F_{\Gamma},\eta_{\Gamma})\rightarrow (B,\delta)$ is a $\Gamma$-equivariant free presentation of $(B,\delta)$, then there is an exact sequence

$0\rightarrow U_{\Gamma}\rightarrow H^{\mathfrak{P}_{\Gamma-e}}_{2}(B,\delta)\rightarrow H^{\mathfrak{P}_{\Gamma-e}}_{2}(G,\mu)\overset{\rho}\rightarrow A/[B,A]_{\Gamma}\rightarrow H^{\mathfrak{P}_{\Gamma-e}}_{1}(B,\delta)\rightarrow H^{\mathfrak{P}_{\Gamma-e}}_{1}(G,\mu)\rightarrow 0$,

where  $U_{\Gamma}$ is the kernel of  $[F_{\Gamma},S_{\Gamma}]_{\Gamma}/[F_{\Gamma},R_{\Gamma}]_{\Gamma}\rightarrow [B,A]_{\Gamma}$, $R_{\Gamma} = Ker \tau_{\Gamma} $ and $S = Ker \theta \tau_{\Gamma}$.
\end{thm}

$Proof$. 1) Let $(F_{*}(G),\eta_{*})\rightarrow (G,\mu)$ and $(F_{*}(B),\beta_{*})\rightarrow (B,\delta)$ be  $\mathfrak{F}$-projective resolutions of $(G,\mu)$ and $(B,\delta)$ respectively induced by $\tau$ and $\theta \tau$. The short exact sequence of augmented pseudo-simplicial groups

   $ ([F_{*},F_{*}]_{\Gamma}\rightarrow [G,G]_{\Gamma})\rightarrow (F_{*}\rightarrow G)\rightarrow (F_{*}/[F_{*},F_{*}]_{\Gamma}\rightarrow G/[G,G]_{\Gamma})$

   yields the short exact sequence

   $0\rightarrow H^{\mathfrak{P}_{\Gamma}}_{2}(G,\mu)\rightarrow \pi_{0}([F_{*},F_{*}]_{\Gamma})\rightarrow [G,G]_{\Gamma}\rightarrow 1$.

   To obtain the required exact sequence it suffices to apply the following commutative diagram with exact rows and columns

$$
\xymatrix{
& 0\ar[d] & 1 \ar[d] & 1 \ar[d] & 0 \ar[d]  \\
0\ar[r] & Ker\alpha'\ar[d]\ar[r] & Ker \alpha \ar[d]\ar[r]^{\kappa} &   A \ar[d]\ar[r] & \ar[d]\ar[r]   A/\kappa(Ker\alpha) \ar[d]\ar[r] & 0\\
0\ar[r] & H^{\mathfrak{P}_{\Gamma}}_{2}(B,\delta)\ar[d]^{\alpha'}\ar[r]&\pi_{0}([F^{B}_{*},F^{B}_{*}]_{\Gamma}) \ar[d]^{\alpha}\ar[r] & B \ar[d]^{\gamma}\ar[r]  & H^{\mathfrak{P}_{\Gamma}}_{1}(B,\delta)
\ar[d]\ar[r] & 0\\
0\ar[r] & H^{\mathfrak{P}_{\Gamma}}_{2}(G,\mu)\ar[r]&\pi_{0}([F^{G}_{*},F^{G}_{*}]_{\Gamma})\ar[d]\ar[r] & G \ar[d]\ar[r] & H^{\mathfrak{P}_{\Gamma}}_{1}(G,\mu)\ar[d]\ar[r] & 0\\
& & 1 & 1 &0
}
$$

It is easily checked that $Ker\alpha = [F,S]_{\Gamma}/[F,R]_{\Gamma}$. Therefore $\kappa(Ker\alpha) = [B,A]_{\Gamma}$ and $Ker\kappa$ is isomorphic to $Ker\alpha'$. The connecting homomorphism $\rho$ is defined in a natural way. This completes the proof of the theorem.

 The Hopf formula for the crossed $\Gamma$-module homology and for the $\Gamma$-equivariant crossed $\Gamma$-module homology follows as a result of Theorem 7.8 Namely

\begin{cor}
1) If $\tau: (F, \eta)\rightarrow (G,\mu)$ is a free presentation of $(G,\mu)$ then

$H^{\mathfrak{P}_{\Gamma}}_{2}(G,\mu) \cong R\cap [F,F]_{\Gamma}/[F,R]_{\Gamma},$
 $R = Ker\tau$.

2)If $\tau_{\Gamma}: (F_{\Gamma}, \eta_{\Gamma})\rightarrow (G,\mu)$ is a $\Gamma$-equivariant free presentation of $(G,\mu)$ then

$H^{\mathfrak{P}_{\Gamma-e}}_{2}(G,\mu) \cong R_{\Gamma}\cap [F_{\Gamma},F_{\Gamma}]_{\Gamma}/[F_{\Gamma},R_{\Gamma}]_{\Gamma},$
 $R_{\Gamma} = Ker\tau_{\Gamma}$.
\end{cor}

$Proof$.
By Theorem 7.8 the $\Gamma$-extension of $(G,\mu)$:

$(R,1)\overset{\sigma}\rightarrow (F,\eta)\rightarrow (G,\mu)$

induces the exact sequence

$H^{\mathfrak{P}_{\Gamma}}_{2}(F,\eta)\rightarrow H^{\mathfrak{P}_{\Gamma}}_{2}(G,\mu)\rightarrow R/[F,R]_{\Gamma}\overset{\sigma'}\rightarrow F/[F,F]_{\Gamma}\rightarrow G/[G,G]_{\Gamma}\rightarrow 0$,

where $\sigma'$ is induced by $\sigma$ and $H^{\mathfrak{P}_{\Gamma}}_{2}(F,\eta) = 0$. since $(F,\eta)$ is a free crossed $\Gamma$-module. Finally we obtain $H^{\mathfrak{P}_{\Gamma}}_{2}(G,\mu)\cong Ker\sigma' = R\cap [F,F]_{\Gamma}/[F,R]_{\Gamma}$.

 For the case 2) of Theorem 7.8 and Corollary 7.9 the proof is completely similar and it is omitted.

\begin{rem}
  Our  results on $\Gamma$-extensions and (co)homology of crossed $\Gamma$-modules can be viewed as a generalization of the homological theory of relative extensions of group epimorphisms of Loday [33] and they are also closely related  to (co)homology of crossed modules of Carrasco, Cegarra and Grandjean, where the category of all crossed modules is considered and the free cotriple (co)homology of crossed modules is used [7].
 \end{rem}

\;

\section{Applications to algebraic K-theory, Galois theory of commutative rings and cohomological dimension of groups}

The first application deals with the connection of the $\Gamma$-equivariant homology of groups and the homology of crossed $\Gamma$-modules with the relative algebraic K-functor $K_{2}(f)$, where $f:\Lambda\rightarrow \Lambda'$ is a surjective homomorphism of rings with unit. For this purpose we recall the definition of $K_{2}(f)$.

\begin{defn}

The relative Steinberg group $St(f)$ of the surjective homomorphism f is the quotient of the free group $F(St(\Lambda)\times Y)$ by the minimal $St(\Lambda)$-equivarient normal subgroup satisfying the relations

$(A_{1})$                          $y^{u}_{ij} y^{v}_{ij} = y^{u+v}_{ij},$

$(B_{1})$                            $x^{\lambda}_{ij}y^{v} = y^{v}_{ij},$

$(B_{2})$                           $x^{\lambda}_{ij}y^{v}_{kl} = y^{v}_{kl}, j\neq k, i\neq l,$

$(B_{3})$                            $x^{\lambda}_{ij}y^{v}_{jk} = y^{\lambda v}_{ik}y^{v}_{jk}, i\neq k,$

$(B_{3'})$                           $x^{\lambda}_{ij}y^{v}_{ki} = y^{-v\lambda}_{kj}y^{v}_{ki}, j\neq k,$

$(C)$                               $x^{\lambda}_{ij}\cdot t = y^{v}_{ij}ty^{-u}_{ij}, t\in F(St(\Lambda)\times Y),$

where Y is the set of $\{y^{u}_{ij}\}$, i,j are positive integers and u belongs to the kernel I of the homomorphism f [33].

\end{defn}

 The homomorphism $\varphi_{f}: St(f)\rightarrow E(\Lambda,I)\subset GL(I)$ is defined by $\varphi_{f}(x\cdot y^{u}_{ij}) = \varphi_{\Lambda}(x)e^{u}_{ij}\varphi_{\Lambda}(x)^{-1}$, where the homomorphism $\varphi_{\Lambda}:ST(\Lambda)\rightarrow E(\Lambda)$ is sending the generator $x^{\lambda}_{ij}$ of $St(\Lambda)$ to $e^{\lambda}_{ij}$. The group $E(\Lambda,I)$ is the direct limit of $\{ E_{n}(\Lambda,I)\},n\rightarrow \infty$, and $E_{n}(\Lambda,I)$ is the normal subgroup of $E_{n}(\Lambda)$ of elementary n-matrices generated by I-elementary matrices of the form $I_{n}+ ve_{ij}$, $v\in I$ and $i\neq j$. The group $E(\Lambda)$ is acting on $E(\Lambda,I)$ by conjugation and it is well known that $E(\Lambda,I)$ is $E(\Lambda)$-perfect.

\begin{defn}

[33] $K_{2}(f) =  Ker \varphi_{f}$ and $Coker \varphi_{f} = K_{1}(f)$.

\end{defn}

The groups $K_{2}(f)$ and $K_{1}(f)$ are also noted $K_{2}(\Lambda, I)$and $K_{1}(\Lambda,I)$ respectively [37,41,25].

Denote D the fiber product $\Lambda \times_{\Lambda'}\Lambda$ with projections $p_{1}:D\rightarrow \Lambda$ and $p_{2}:D\rightarrow \Lambda$. Let St(I) be the kernel of $St(p_{1})$ and let $C(I)= [St(p_{1}),St(p_{2})]$. There is a homomorphism $\mu: St(I)/C(I)\rightarrow St(\Lambda)$ induced by $St(p_{2})$ on St(I). In [33] it is shown that the set of relations defining the group St(I) given by Swan [46] is equivalent to the set of relations $(A_{1},B_{1},B_{2},B_{3},B_{3'})$ implying the isomorphisms $\theta: St(f)\rightarrow St(I)/C(I)$, $K_{2}(f)\overset{\cong}\rightarrow K_{2}(I)/C(I)$, and the sequence

\begin{equation}
0\rightarrow Ker(\mu\theta)\rightarrow St(f)\overset{\mu\theta}\rightarrow St(\Lambda)\overset{St(f)}\rightarrow St(\Lambda')\rightarrow 1
\end{equation}
is the universal relative extension of $(St(\Lambda'),St(\Lambda))$.

According to results of [33] the following short exact sequence is provided

\begin{equation}
0\rightarrow K_{2}(f)\rightarrow St(f)\overset{\varphi_{f}}\rightarrow E(\Lambda,I)\rightarrow 1,
\end{equation}

where $\varphi_{f} = \varphi_{\Lambda}\mu \vartheta$ (see also [30]).

\begin{thm}

There is an exact sequence

$0\rightarrow [P,S]_{St(\Lambda)}/[P,R]_{St(\Lambda)}\rightarrow H^{\mathfrak{P}_{St(\lambda)}}_{2}(St(f)/St(\Lambda)(\varphi_{f}))\rightarrow H^{\mathfrak{P}_{St(\lambda)}}_{2}(E(\Lambda, I))\rightarrow$

$\rightarrow K_{2}(f)/St(\Lambda)(\varphi_{f}))\rightarrow 0,$

\end{thm}

where  $\alpha: P\rightarrow St(f)$ is a $St(\Lambda)$-projective presentation of St(f),$R = Ker \alpha$ and $S = Ker \varphi_{f}\alpha.$

 $Proof$. Consider the normal subgroup of St(f) generated by the elements $^{\gamma}x\cdot x^{-1}, x\in St(f), \gamma\in St(\Lambda)$, such that $\varphi_{f}(^{\gamma}x\cdot x^{-1}) = 1.$. This subgroup is denoted $St(\Lambda)(\varphi_{f})$. By Corollary 3.3 this yields the exact sequence

 $0\rightarrow K_{2}(f)/St(\Lambda)(\varphi_{f}))\rightarrow St(f)/St(\Lambda)(\varphi_{f}))\rightarrow E(\Lambda,I)\rightarrow 1$

which is a central $St(\Lambda)$-equivariant extension of $E(\Lambda,I)$ having a $St(\Lambda)$-section map. The group St(f)is $St(\Lambda)$-perfect [33] implying $St(f)/St(\Lambda)(\varphi_{f}))$ is also $St(\Lambda)$-perfect and therefore $H^{\mathfrak{P}_{St(\lambda)}}_{1}(St(f)/St(\Lambda)(\varphi_{f}))) = 0$. It remains to apply Theorem 2.8 to get the required exact sequence. This completes the proof.

\begin{rem}
The exact sequence of Theorem 8.3 can be replaced by the exact sequence

$0\rightarrow [P,S]_{St(\Lambda)}/[P,R]_{St(\Lambda)}\rightarrow H^{\mathfrak{P}_{St(\lambda)}}_{2}(St(f)/St(\Lambda)(\varphi_{f}))\rightarrow \tilde{H}^{\mathfrak{P}_{St(\lambda)}}_{2}(E(\Lambda, I))\rightarrow K_{2}(f)\rightarrow 0$,

 where $\tilde{H}^{\mathfrak{P}_{St(\lambda)}}_{2}(E(\Lambda, I))$ denotes the fiber product

 $H^{\mathfrak{P}_{St(\lambda)}}_{2}(E(\Lambda, I))\times_{K_{2}(f)/St(\Lambda)(\varphi_{f})} K_{2}(f)$.
\end{rem}

We are now going to establish the relation of the homology of crossed $\Gamma$-modules with the relative algebraic K-functor $K_{2}(f)$.

The short exact sequence (8.2) induces the following $E(\Lambda)$-extension of the inclusion crossed $E(\Lambda)$-module $(E(\Lambda,I), \sigma)$, $\sigma:E(\Lambda,I)\hookrightarrow E(\Lambda)$:

$0\rightarrow (K_{2}(f),1)\rightarrow (St(f),\varphi_{f})\rightarrow (E(\Lambda,I),\sigma)\rightarrow 1.$

Take the quotient St'(f) of St(f) by the normal subgroup generated by the elements $^{\gamma}x\cdot x^{-1}, \gamma\in Ker \varphi_{\Lambda}, x\in St(f),$ implying the short exact sequence

$0\rightarrow K_{2}(f)\rightarrow St'(f)\overset{\varphi'_{f}}\rightarrow E(\Lambda,I)\rightarrow 1$

of $E(\Lambda)$-modules, $E(\Lambda)$ is trivially acting on $K_{2}(f)$ and its action on St'(f) is realized via the homomorphism $\varphi_{\Lambda}$.

Finally we obtain a central $E(\Lambda)$-extension of the inclusion crossed $E(\Lambda)$-module $(E(\Lambda,I),\sigma)$, $\sigma: E(\Lambda,I)\hookrightarrow E(\Lambda)$,
\begin{equation}
0\rightarrow (K_{2}(f),1)\rightarrow (St'(f),\varphi'_{f})\rightarrow (E(\Lambda,I),\sigma)\rightarrow 1,
\end{equation}
where $\varphi'_{f}$ is induced by $\varphi_{\Lambda}\mu \theta$.

The sequence

\begin{equation}
0\rightarrow (K_{2}(f)/E(\Lambda) (\varphi'_{f}),1)\rightarrow (St'(f)/E(\Lambda) (\varphi'_{f}),\sigma \varphi'_{f})\rightarrow (E(\Lambda,I),\sigma)\rightarrow 1
\end{equation}

is $E(\Lambda)$-equivariant extension of the inclusion crossed $E(\Lambda)$-module $(E(\Lambda,I),\sigma)$. It is evident that the crossed $E(\Lambda)$-module $(St'(f), \sigma \varphi'_{f})$ is $E(\Lambda)$-perfect and therefore the crossed $(E(\Lambda)$-module $(St'(f)/E(\Lambda)(\varphi'_{f}), \sigma \varphi'_{f})$ is also $E(\Lambda)$-perfect implying $H^{\mathfrak{P}_{E(\Lambda)-e}}_{1}(St'(f)/E(\Lambda)(\varphi'_{f}), \sigma \varphi'_{f}))
 = 0.$

\begin{thm}
1) The sequence of $St(\Lambda)-modules$
$$
0\rightarrow (Ker(\mu\theta),1)\rightarrow (St(f),\mu\theta)\rightarrow (Ker St(f),\sigma)\rightarrow 1
$$

is the universal $St(\Lambda$)-extension of the inclusion crossed $St(\Lambda)$-module $(Ker St(f),\sigma)$ and there is an isomorphism $H^{\mathfrak{P}_{St(\lambda)}}_{2}(Ker St(f),\sigma)\cong Ker(\mu\theta)$.

2) The sequence (8.4) is the $E(\Lambda)$-equivariant extension of the inclusion crossed $E(\Lambda)$-module $(E(\Lambda,I),\sigma)$ and there is an isomorphism

$ H^{\mathfrak{P}_{E(\Lambda)-e}}_{2}(E(\Lambda,I),\sigma) \overset{\cong}\rightarrow K_{2}(f)/E(\Lambda) (\varphi'_{f}).$

\end{thm}

$Proof$. 1) As noted above it is proven in [33] that the group St(f) is $St(\Lambda)$-perfect. Therefore the crossed $St(\Lambda)$-module $(St(f,\mu\theta)$ is also $St(\Lambda)$-perfect and $H^{\mathfrak{P}_{St(\lambda)}}_{1}(St(f),\mu\theta) = 0$. Since the sequence (8.1) is the universal relative extension of $(St(\Lambda'), St(\Lambda))$ it follows that every central $St(\Lambda)$-extension of the crossed $St(\Lambda)$-module $(St(f), \mu\theta)$ splits implying $H^{\mathfrak{P}_{St(\lambda)}}_{2}(St(f),\mu\theta) = 0$. It remains to apply the first part of Theorem 7.8 to get the required isomorphism.

2) It is evident that the crossed $E(\Lambda)$-module $(St'(f),\varphi'_{f})$ is $E(\Lambda)$-perfect. Therefore the crossed $E(\Lambda)$-module $(St'(f)/E(\Lambda) (\varphi'_{f}),\sigma \varphi'_{f})$ is also $E(\Lambda)$-perfect and $H^{\mathfrak{P}_{E(\lambda)-e}}_{1}(St'(f)/E(\Lambda) (\varphi'_{f}),\sigma \varphi'_{f}) = 0$.

Let

$$
0\rightarrow (A,1)\rightarrow (U',\eta')\overset{\beta'}\rightarrow (St'(f)/E(\Lambda) (\varphi'_{f}),\sigma \varphi'_{f})\rightarrow 1
$$

be a central $E(\Lambda)$-equivariant extension of $(St'(f)/E(\Lambda) (\varphi'_{f}),\sigma \varphi'_{f})$. Consider the fiber product

$$
\xymatrix{
 D' \ar[d]^{q'_{2}}\ar[r]^{q'_{1}} & (St(f'),\varphi'_{f})\ar[d]^{g'_{1}}\\
(U',\eta')\ar[r]^{\beta'}& St'(f)/E(\Lambda) (\varphi'_{f}),\sigma \varphi'_{f}).\\
 }
$$,
where $g'_{1}$ is induced by the canonical homomorphism $St(f')\rightarrow St(f')/E(\Lambda)(\varphi'_{f})$.
The sequence
$0\rightarrow (A,1)\rightarrow D'\overset{q'_{1}}\rightarrow (St(f'),\varphi'_{f})\rightarrow 1$ is $E(\Lambda)$-equivariant extension of $(St'(f),\varphi'_{f})$ which becomes a $St(\Lambda)$-equivariant extension of $(St'(f),\varphi'_{f})$ via the homomorphism $\varphi_{\Lambda}:St(\Lambda)\rightarrow E(\Lambda)$.

Now take the fiber product

$$
\xymatrix{
 D \ar[d]^{q_{2}}\ar[r]^{q_{1}} & (St(f),\varphi_{f})\ar[d]^{g_{1}}\\
D'\ar[r]^{q'_{1}}& (St(f'),\varphi'_{f}).\\
 }
$$

The $St(\Lambda)$-equivariant extension $D\overset{q_{1}}\rightarrow (St(f),\varphi_{f})$ is $St(\Lambda)$-splitting and therefore it $E(\Lambda)$-splits too. Let $\sigma_{1}: (St(f),\varphi_{f})\rightarrow D$ be the splitting homomorphism implying the homomorphism $q_{2}\sigma_{1}:(St(f),\mu\theta)\rightarrow D'$
of crossed $E(\Lambda)$-modules such that $q'_{1} q_{2}\sigma_{1}:(St(f),\mu\theta) = g_{1}$. For $^{\gamma}x\cdot x^{-1} \in Ker\varphi_{\Lambda}$ we have $g_{1}(^{\gamma}x\cdot x^{-1}) = 1$ implying the equality $q'_{1}(^{\gamma}q_{2}\sigma_{1}(x)\cdot q_{2}\sigma_{1}(x)^{-1}) = 1.$ By Theorem 3.2 the $E(\Lambda)$-equivariant extension $D'\overset{q'_{1}}\rightarrow (St(f'), \varphi'_{f})$ has the E-property implying the equalities $q_{2}\sigma_{1}(^{\gamma}x\cdot x^{-1}) = ^{\gamma}q_{2}\sigma_{1}(x)\cdot q_{2}\sigma_{1}(x)^{-1} = 1.$ Therefore the homomorphism $q_{2}\sigma_{1}$ sends to the unit the normal subgroup of St(f) generated by the elements $^{\gamma}x\cdot x^{-1} \in Ker\varphi_{\Lambda}$, inducing the $E(\Lambda)$-homomorphism $\sigma'_{1}:(St'(f),\varphi'_{f})\rightarrow D'$ such that $q'_{1}\sigma'_{1} = 1.$

Starting with the first diagram and with the splitting homomorphism $\sigma'_{1}$ it is easily shown by the same line of argumentation as for the previous case that the $E(\Lambda)$-equivariant extension $\beta':(U',\eta')\rightarrow (St'(f)/E(\Lambda)(\varphi'_{f}),\sigma \varphi'_{f})$ splits. It follows that $H^{\mathfrak{P}_{E(\lambda)-e}}_{2}(St'(f)/E(\Lambda) (\varphi'_{f}),\sigma \varphi'_{f}) = 0$ and it remains to apply the second part of Theorem 7.8. This completes the proof of the theorem.

The second application concerns the investigation of the relationship between the equivariant symbol group of non commutative local rings and the Milnor algebraic K-functor $K_{2}$ by using the $\Gamma$-homology of groups and the homology of crossed $\Gamma$-modules that will extend the well known Matsumoto's theorem for fields [34].

Let A be a unital ring and Sym(A) be the symbol group of the ring A generated by the elements $\{u,v\}$, $u,v \in A^{*}$, satisfying the following relations

$(S_{0})$     $\{u,1-u\} = 1, u\neq 1, 1-u \in A^{*},$

$(S_{1})$     $\{uu',v\} = \{u,v\}\{u',v\}$,

$(S_{2})$     $\{u,vv'\} = \{u,v\}\{u,v'\}$.
where $A^{*}$ denotes the multiplicative group of invertible elements of the ring A [ ]. By Matsumoto's theorem the groups Sym(A) and $K_{2}(A)$ are isomorphic when A is a field [34].

For our purpose it is necessary to introduce the notion of equivariant symbol group.

\begin{defn}
  For a ring A with unit the equivariant symbol group $Sym^{A^{*}}(A)$ is defined as the group $(Sym(A))_{A^{*}}$.
\end{defn}

The symbol group Sym(A) becomes $A^{*}$-group by the action $^{u}\{v,w\} = \{uvu^{-1},uwu^{-1}\}$. Therefore the equivariant symbol group of an unital commutative ring coincides with its symbol group.

Now assume the ring A is a non commutative local ring such that $A/Rad(A) \neq F_{2}$. Consider the group U(A) generated by the elements $\langle u,v\rangle$,$u,v \in A^{*}$,$A^{*}$ , satisfying the following relations

$(U_{0})$        $\langle u, 1-u\rangle = 1, u\neq 1, 1-u \in A^{*},$

$(U_{1})$        $\langle uv,w \rangle = ^{u}\langle v,w \rangle \langle u,w\rangle$,

$(U_{2})$        $\langle u,vw\rangle \langle v,wu\rangle \langle w,uv\rangle = 1$,

where $^{u}\langle v,w \rangle = \langle uvu^{-1},uwu^{-1} \rangle$ [21].

The group U(A) becomes $A^{*}$-group with respect to this action and results of [21] show us that there is a surjective $A^{*}$-homomorphism $U(A)\rightarrow Sym(A)$. In addition there is a short exact sequence of $A^{*}$-groups relating U(A) with $K_{2}(A)$ [21]:

\begin{equation}
0\rightarrow K_{2}(A)\rightarrow U(A)\overset{\tau}\rightarrow [A^{*},A^{*}]\rightarrow 1,
\end{equation}
where $A^{*}$ acts trivially on $K_{2}(A)$ and by conjugation on $[A^{*},A^{*}]$, and $\tau(\langle u,v\rangle) = [u,v]$. Moreover the sequence (8.5) induces a central $A^{*}$-extension of the inclusion crossed $A^{*}$-module $([A^{*},A^{*}],i)$:

\begin{equation}
0\rightarrow (K_{2}(A),1)\rightarrow (U(A),\tau)\overset{i\tau}\rightarrow ([A^{*},A^{*}],i)\rightarrow 1.
\end{equation}

We will need the corresponding $A^{*}$-equivariant versions of these two sequences. Namely,

\begin{equation}
 0\rightarrow (K_{2}(A))/A^{*}(\tau)\rightarrow U'(A)\overset{\tau'}\rightarrow [A^{*},A^{*}]\rightarrow 1,
 \end{equation}
 \begin{equation}
 0\rightarrow (K_{2}(A))/A^{*}(\tau),1)\rightarrow (U'(A),i\tau')\overset{\tau'}\rightarrow ([A^{*},A^{*}],i)\rightarrow 1,
\end{equation}
where $U'(A) = (U(A))/A^{*}(\tau)$ and $\tau'$ is induced by $\tau$.

The subgroup $A^{*}(\tau)$ of $K_{2}(A)$ is generated by the elements $\langle \gamma, x\rangle$ such that $\gamma x = x\gamma$, where $\gamma\in A^{*}, x\in [A^{*},A^{*}]$. In effect, let $\prod \langle u_{i},v_{i}\rangle$ be an element of U(A). One has $\prod ^{\gamma}\langle u_{i},v_{i}\rangle = \prod \langle \gamma,[u_{i},v_{i}]\rangle \langle u_{i},v_{i}\rangle = \langle\gamma,\prod[u_{i},v_{i}]\rangle \prod\langle u_{i},v_{i}\rangle$. Thus $^{\gamma}(\prod\langle u_{i},v_{i}\rangle)\cdot (\prod\langle u_{i},v_{i}\rangle)^{-1} = \langle \gamma,\prod [u_{i},v_{i}]\rangle$ and  the equality $\tau( ^{\gamma}(\prod\langle u_{i},v_{i}\rangle)\cdot (\prod\langle u_{i},v_{i}\rangle)^{-1}) = 1$ implies $\tau(^{\gamma}(\prod\langle u_{i},v_{i}\rangle)\cdot (\prod\langle u_{i},v_{i}\rangle)^{-1}) = [^{\gamma}\prod \langle u_{i},v_{i}] = 1$.

By this way the defining relations for the group U'(A) have been also provided as follows:

$(U_{0})$        $\langle u, 1-u\rangle = 1, u\neq 1, 1-u \in A^{*},$

$(U_{1})$        $\langle uv,w \rangle = ^{u}\langle v,w \rangle \langle u,w\rangle$,

$(U_{2})$        $\langle u,vw\rangle \langle v,wu\rangle \langle w,uv\rangle = 1$,

$(U_{3})$        $uv = vu, u\in A^{*}, v\in [A^{*},A^{*}]$.

\begin{thm}
There is an exact sequence
$0\rightarrow [P,S]_{A^{*}}/[P,R]_{A^{*}}\rightarrow H^{A^{*}}_{2}(U'(A))\rightarrow H^{A^{*}}_{2}([A^{*},A^{*}])\rightarrow K_{2}(A)/A^{*}(\tau)\rightarrow Sym^{A^{*}}(A)\rightarrow [A^{*},A^{*}]/[A^{*},A^{*}]_{A^{*}}\rightarrow 0$,
\end{thm}
where $\alpha:P\rightarrow U'(A)$ is $A^{*}$-projective presentation of U'(A), $R = Ker\alpha$ and $S = Ker\tau'\alpha$.

$Proof$.  First of all it will be proved that the groups $H_{1}^{A^{*}}(U'(A))$ and $Sym^{A^{*}}(A)$ are isomorphic. Consider the system of relations $(S_{0},S_{1},S'_{2})$ which is equivalent to the system $(S_{0},S_{1},S_{2})$ of defining relations for the symbol group Sym(A)[21]. It is easily checked that the system of relations $((S_{0},S_{1},S'_{2},L_{4},L_{5})$ is equivalent to the system $(U_{0},U_{1},U_{2},M_{4},M_{5})$, where $(L_{4}) : ^{w}\{u,v\} = \{u,v\}, (L_{5}): \{u,v\}\{u',v'\} = \{u',v'\}\{u,v\}$,and $(M_{4}) : ^{w}\langle u,v\rangle = \langle u,v\rangle, (M_{5}): \langle u,v\rangle \langle u',v'\rangle = \langle u',v'\rangle \langle u,v\rangle$. For the group U(A) the relation $\langle u,v\rangle\langle u',v'\rangle = ^{[u,v]}\langle u',v'\rangle\langle u,v\rangle$ holds [21], implying the group $U(A)_{A^{*}}$ defined by the system of relations $((U_{0},U_{1},U_{2},M_{4})$ is abelian. Therefore the group $Sym(A)_{A^{*}}$ is also abelian and this yields the following sequence of equivalences of defining systems:

$(S_{0},S_{1},S_{2},L_{4}) \approx (S_{0},S_{1},S'_{2},L_{4}) \approx ((S_{0},S_{1},S'_{2},L_{4},L_{5}) \approx (U_{0},U_{1},U_{2},M_{4},M_{5}) \approx  (U_{0},U_{1},U_{2},M_{4})$

showing the isomorphism of groups $Sym(A)_{A^{*}}$ and $U(A)_{A^{*}}$. Finally this induces the following isomorphisms $H_{1}^{A^{*}}(U'(A)) \cong U'(A)_{A^{*}} \cong U(A)_{A^{*}} \cong Sym(A)_{A^{*}}$. It remains to apply Theorem 2.8 for the sequence (8.7). This completes the proof of the theorem.

\begin{cor}
(1) As a consequence of this theorem there is an exact sequence

$0\rightarrow [P,S]_{A^{*}}/[P,R]_{A^{*}}\rightarrow H^{A^{*}}_{2}(U'(A))\rightarrow \tilde{H}^{A^{*}}_{2}([A^{*},A^{*}])$

$\rightarrow K_{2}(A)\rightarrow  Sym^{A^{*}}(A)\rightarrow [A^{*},A^{*}]/[A^{*},A^{*}]_{A^{*}}\rightarrow 0,$

and

(2) if $[A^{*},A^{*}]$ is quasi-perfect then the sequence

$0\rightarrow [P,S]_{A^{*}}/[P,R]_{A^{*}}\rightarrow H^{A^{*}}_{2}(U'(A))\rightarrow \tilde{H}^{A^{*}}_{2}([A^{*},A^{*}])$

$\rightarrow K_{2}(A)\rightarrow  Sym^{A^{*}}(A)\rightarrow 0$

is exact.
\end{cor}
where $\tilde{H}^{A^{*}}_{2}([A^{*},A^{*}])$ is the fiber product

 $H^{A^{*}}_{2}([A^{*},A^{*})])\times_{K_{2}(A)/A^{*}(\tau)} K_{2}(A)$.

The sequence of Corollary 8.8,(2) generalizes the exact sequence

$H_{2}(U(D))\rightarrow H_{2}([D^{*},D^{*}])\rightarrow K_{2}(D)\rightarrow Sym(D)\rightarrow 1$

given in [2] for a non commutative divisor ring D such that $[D^{*},D^{*}]$ is perfect, and the sequence of Corollary 8.8(1) can be considered as an abelian version of the exact sequence of Guin [21]

$(A^{*})^{ab}\otimes_{Z}K_{2}(A)\rightarrow \bar{H}_{1}(A^{*},U(A))\rightarrow \bar{H}_{1}(A^{*},[A^{*},A^{*}])\rightarrow K_{2}(A)$

$\rightarrow Sym(A)\rightarrow [A^{*},A^{*}]/[A^{*},[A^{*},A^{*}]]\rightarrow 1$,

where $\bar{H}_{1}$ is the first non abelian homomolgy of groups with coefficients in crossed modules.

\begin{rem}
Guin's low dimensional non-abelian group homoloy with coefficients in crossed modules is closely related to integral homology of crossed modules.
Let $(A,\delta)$ be a crossed G-module and $\bar{H_{0}}(G,A), \bar{H_{1}}(G,A)$ denote Guin's group homology with coefficients in the crossed $G$-module $(A,\delta)$. In [5] it is shown that there is a group homomorphism $\varphi: G\bigotimes A\rightarrow A$, $\varphi(g\otimes a) =  ^{g}a\cdot a^{-1}$, where $G\bigotimes A$ is the non-abelian tensor product of Brown-Loday, and $(G\bigotimes A, \varphi)$ is a crossed A-module. Then $\bar{H_{0}}(G,A) = coker \varphi$ and $\bar{H_{1}}(G,A) = Ker \varphi$. It is easily seen that there is an isomorphism of the abelianization $\bar{H}^{ab}_{0}(G,A)$ with $H^{\mathfrak{P}}_{1}(A,\delta)$ and there is an exact sequence

 $ H^{\mathfrak{P}}_{2}(G\bigotimes A,\varphi)\rightarrow H^{\mathfrak{P}}_{2}(\Gamma A,\sigma)\rightarrow \bar{H_{1}}(G,A)\rightarrow G\bigotimes A/([G\bigotimes A,G\bigotimes A] )_{A}\rightarrow \Gamma A/([\Gamma A,\Gamma A])_{A}\rightarrow 0$

 induced by the exact central sequence of crossed A-modules

 $0\rightarrow (\bar{H_{1}}(G,A),1)\rightarrow (G\bigotimes A, \varphi)\rightarrow (\Gamma A,\sigma)\rightarrow 1$.
\end{rem}

\ {The next application concerns the connection of $\Gamma$-equivariant derived functors and $\Gamma$-equivariant cohomology of groups with Galois theory of commutative rings.

Let R be a commutative ring with no nontrivial idempotents in which the prime number p is invertible. Let $C_{p^{n}}$ be the cyclic group of order $p^{n}$ and $S_{n}$ be the splitting ring of the polynomial $x^{p^{n}}-1$ of R[x]. The following notation is also introduced: $S^{*}_{n}$ is the group of invertible elements and $\mu_{p^{n}}$ is the group of n-roots of 1 in the splitting ring $S_{n}$, $NB(R,C_{p^{n}})$ is the set of isomorphism classes of Galois extensions with normal basis of R with Galois group $C_{p^{n}}$ and $\Gamma_{n}$ is the Galois group of $S_{n}$.

\begin{thm}
There are bijections

1. $NB(R,C_{p^{n}}) \cong$  $Ext^{1}_{\mathds{Z},\Gamma_{n}}(\mu_{p^{n}},S^{*}_{n})$, where $\Gamma_{n}$ is trivially acting on $\mathds{Z}$.

2.  $NB(R,C_{p^{n}}) \cong$ $H^{2}_{\Gamma_{n}}(\mu_{p^{n}},S^{*}_{n})$, where $S^{*}_{n}$ is a trivial $\mu_{p^{n}}$-module.
\end{thm}

$Proof$. First of all it is necessary to show that there is a bijection of $Ext^{1}_{\Lambda,\Gamma}(L,M)$ with $E^{1,\Lambda}_{\Gamma}(L,M)$ which is the set of isomorphism classes of short exact sequences $0\rightarrow M\rightarrow X\rightarrow L\rightarrow 0$ of $\Gamma$-equivariant $\Lambda$-modules having a $\Gamma$-section map.

Consider the standard $\Gamma$-projective resolution of the $\Gamma$-equivariant $\Lambda$-module L:

$ ... \rightarrow F_{n}\overset{\alpha_{n}}\rightarrow F_{n-1}\overset{\alpha_{n-1}}\rightarrow ... \overset{\alpha_{2}}\rightarrow F_{1}\overset{\alpha_{1}} \rightarrow F_{0}(L)\overset{\tau}\rightarrow L\rightarrow 0$

where $F_{0}(L)$ is the free $\Lambda$-module with basis $|l|$, $l\in L$, $\tau$ is the canonical surjective homomorphism and the action of $\Gamma$ on $F_{0}(L)$ is given by $^{\gamma}|l| = |^{\gamma}l|$, $l\in L$, $F_{1} = F_{0}(Ker \tau)$ and by induction on n the relatively free $\Gamma$-equivariant $\Lambda$-module $F_{n}$ , $n> 1,$ is defined as $F_{0}(Ker \alpha_{n-1})$, the homomorphism $\alpha_{n}$ is induced by the canonical homomorphism $F_{0}(Ker \alpha_{n-1})\rightarrow Ker \alpha_{n-1}$. Then $Ext^{n}_{\Lambda,\Gamma}(L,M)$ , $n\geq 0$, is the n-th homology group of the chain complex

$0\rightarrow Hom^{\Gamma}_{\Lambda}(F_{0}(L),M)\rightarrow  Hom^{\Gamma}_{\Lambda}(F_{1},M)\rightarrow ... \rightarrow Hom^{\Gamma}_{\Lambda}(F_{n-1},M)\rightarrow Hom^{\Gamma}_{\Lambda}(F_{n},M)\rightarrow ...$.

It is evident that there is an ismorphism $Ext^{0}_{\Lambda,\Gamma}(L,M) \cong Hom^{\Gamma}_{\Lambda}(L,M)$. Let $E: 0\rightarrow M\rightarrow X\rightarrow L\rightarrow 0$ be a short exact sequence of $\Gamma$-equivariant $\Lambda$-modules having $\Gamma$-section map $\sigma:L\rightarrow X$, $[E]\in E^{1,\Lambda}_{\Gamma}(L,M)$. The homomorphism $\tau$ induces a $\Lambda$-homorphism $F_{0}(L)\rightarrow X$ compatible with the action of $\Gamma$ sending $|l|$ to $\sigma(l)$. Let f denote its restriction on $Ker \tau$ and let $[f\alpha_{0}]$ be the class of the homomorphism $f\alpha_{0}:F_{1}\rightarrow M$. Therefore we obtain a map from $E^{1,\Lambda}_{\Gamma}(L,M)$ to $Ext^{1}_{\Lambda,\Gamma}(L,M)$ sending [E] to $[f\alpha_{0}]$. Conversely, if $[g]\in Ext^{1}_{\Lambda,\Gamma}(L,M)$, then g induces a $\Lambda$-homomorphism $g':Ker \tau\rightarrow M$ compatible with the action of $\Gamma$. Similarly to the classical case (when $\Gamma$ is acting trivially on $\Lambda$) by using the homomorphism g' we can construct a short exact sequence of $\Gamma$-equivariant $\Lambda$-modules with $\Gamma$-section map as follows. Take the sum $M\oplus F_{0}(L)$ which is a $\Gamma$-equivariant $\Lambda$-module with componentwise action of $\Gamma$ and its quotient Y by the submodule generated by the elements $(-g'(x), \vartheta(x))$, $x\in Ker \tau$, where $\vartheta:Ker \tau\rightarrow F_{0}(L)$ is the inclusion map. This yields the needed short exact sequence $E_{g}: 0\rightarrow M \overset{\beta}\rightarrow Y\overset{\eta}\rightarrow L\rightarrow 0$ of $\Gamma$-equivariant $\Lambda$-modules, where $\beta(m) = [(m,0)]$, $\eta([m,x]) = \tau(x)$, and the $\Gamma$-section map is $\delta(l) = [(0,|l|)]$. Thus we obtain a map $Ext^{1}_{\Lambda,\Gamma}(L,M)\rightarrow E^{1,\Lambda}_{\Gamma}(L,M)$ sending [g] to $[E_{g}]$. It is easily checked that both maps induced by $[E]\rightarrowtail [f\alpha_{0}]$ and $[g]\rightarrowtail [E_{g}]$ respectively are inverse to each other.

In [29] there is the following formula of G.Janelidze

$NB(R,C_{p^{n}}) = Ext^{1}_{\underline{S_{n}}}(Hom(J,U_{n}(R))),U(R_{n}))$, where $J = C_{p^{n}}$, $U(R_{n}) = S^{*}_{n}$, $U_{n}(R) = \mu_{p^[n]}$ and $\underline{S_{n}}$ is the category of $\Gamma_{n}$-sets.

As noted by Greither [20] this beautiful formula allowed us to establish the bijection of $RNB(R,C_{p^{n}})$ with the set of isomorphism classes of short exact sequences $0\rightarrow S^{*}_{n}\rightarrow X\rightarrow \mu_{p^{n}}\rightarrow 0$ of $\Gamma_{n}$-modules having $\Gamma_{n}$-section map. This completes the proof of the first bijection.

For the second bijection it suffices to remark that suppose $0\rightarrow M\rightarrow X\rightarrow L\rightarrow 0$  is a short sequence of $\Gamma$-groups with $\Gamma$-section map, where L is abelian and M is a $\Gamma$-equivariant L-module. If L is a cyclic group trivially acting on M and $\Gamma$ is trivially acting on $\mathds{Z}$, then the group X is abelian and the considered sequence is a short exact sequence of $\Gamma$-equivariant L-modules with $\Gamma$-section map. In that case by [27, Theorem 20] this implies the bijection $E^{1,\mathds{Z}}_{\Gamma}(L,M) \cong H^{2}_{\Gamma}(L,M)$. This completes the proof of the theorem.

Finally, the relation of $\Gamma$-equivariant cohomology of $\Gamma$-groups with equivariant dimensions of groups with operators will be established, particularly with the equivariant cohomological dimension of $\Gamma$-groups.

Recently in [19] the important and well known theorems of Eilenberg-Ganea [14] and Stallings-Swan [43,45] relating the cohomological dimension, the geometric dimension and the Lusternik-Schnirelmann category have been extended to the setting of $\Gamma$-groups as follows:

1) Equivariant Eilenberg - Ganea Theorem:

  Let G be a $\Gamma$-group, where $\Gamma$ is finite. Then the chain of inequalities
$$
cd_{\Gamma}(G) \leq cat_{\Gamma}(G) \leq gd_{\Gamma}(G) \leq sup\{3, cd_{\Gamma}(G)\}
$$

is satisfied. Furthermore, if $cd_{\Gamma}(G) = 2$ then $cat_{\Gamma}(G) = 2$.

2)  Equivariant Stallings - Swan Theorem:

Let G be a $\Gamma$-group, where $\Gamma$ is finite. The following equalities are equivalent:
$$
(1) gd_{\Gamma}(G) = 1,
(2) cat_{\Gamma}(G) = 1,
(3) cd_{\Gamma}(G) = 1,
$$
(4) G is a non-trivial $\Gamma$-free group.

For this purpose the equivariant version of these three quantities have been provided and the equivariant group cohomology has been introduced that is the generalization of the $\Gamma$-equivariant cohomology of $\Gamma$-groups allowing a wider class of coefficients. It is defined as the group cohomology $H^{*}(O_{\mathcal{G}}(G\rtimes \Gamma),M) = Ext^{*}_{O_{\mathcal{G}}(G\rtimes \Gamma)}(\underline{{\mathcal{Z}}},M)$, where $\mathcal{G}$ denotes the family of subgroups of $G\rtimes \Gamma$ which are conjugate to a subgroup of $\Gamma$ and $O_{\mathcal{G}}(G\rtimes \Gamma)$ is the orbit category whose objects are the $\Gamma$-sets $(G\rtimes\Gamma) /H$ for $H\in \mathcal{G}$ and morphisms are $\Gamma$-maps. An $O_{\mathcal{G}}(G\rtimes \Gamma)$-module is a contravariant functor from the category $O_{\mathcal{G}}(G\rtimes \Gamma)$ to the category of abelian groups and $\underline{\mathbb{Z}}$ is the constant functor with value ${\mathbb{Z}}$. The equivariant cohomological dimension $cd_{\Gamma}(G)$ of a $\Gamma$-group G is defined as the least dimension d such that $H^{d+1}(O_{\mathcal{G}}(G\rtimes \Gamma),M) = 0$ for all $O_{\mathcal{G}}(G\rtimes \Gamma)$-modules M.

It would be natural to introduce another algebraic cohomological dimension $cd^{\Gamma}(G)$ of a $\Gamma$-group G based on the $\Gamma$-equivariant cohomology of $\Gamma$-groups as follows:

\begin{defn}
   The cohomological dimension $cd^{\Gamma}(G)$ of the $\Gamma$- group G is the least dimension d such that $H^{d+1}_{\Gamma}(G,M) = 0$ for all $(G\rtimes \Gamma)$-modules M.
 \end{defn}

In [19, Remark 9.1] it is notified that the $\Gamma$-equivariant cohomology of $\Gamma$-groups is the relative group cohomology in the sense of Hochschild [22] and Adamson [1] (see also Benson [4]) enhancing the interest to this equivariant group cohomology. The $\Gamma$-equivariant group cohomology $H^{*}_{\Gamma}(G,M)$ is isomorphic to $H^{*}(G\rtimes \Gamma),\Gamma; M)$ and there is an isomorphism $H^{*}(G\rtimes \Gamma),\Gamma; M) \cong H^{*}(O_{\mathcal{G}}(G\rtimes \Gamma),M^{-})$, where $M^{-}$ is a $O_{\mathcal{G}}(G\rtimes \Gamma)$-module induced by M [39].

Therefore we have the inequality $cd^{\Gamma}(G)\leq cd_{\Gamma}(G)$.

Conjecture:  $cd^{\Gamma}(G)\neq cd_{\Gamma}(G)$.

Problem: Prove the above mentioned  Eilenberg - Ganea and Stallings - Swan theorems in the $\Gamma$-equivariant group cohomology settings involving the cohomological dimension $ch^{\Gamma}(G)$ given in Definition 9.11 and the relevant quantities: the geometric dimension and the Lusternik - Schnirelmann category.

\;



\begin{bibdiv}
\begin{biblist}



\bib{Ad}{article}{
author={Adamson I.T},
title={Cohomology theory for non-normal subgroups and non-normal fields},
journal={Glasgow Math. Ass.},
volume={1954(2)}
date={1954},
pages={66-76},
}

\bib{AlDe}{article}{
author={Alperin R.},
author={Denis K.},
title={$K_{2}$ of quaternion algebra},
journal={J.Algebra},
volume={56},
date={1979},
pages={262-273},
}

\bib{Ba}{article}{
author={Baues H.J.},
author={Thom, A.},
title={Combinatorial homology and 4-dimensional complexes},
journal={De Gruyter, Berlin},
date={1991},

}


\bib{Be}{article}{
author={Benson D.J},
title={Representations and Cohomology I},
journal={Cambridge Studies in Advances Mathematics},
volume={30},
date={1995},
}

\bib{BrLo}{article}{
author={Brown R.},
author={J.-L.Loday},
title={Excision homotopique en basse dimension},
journal={C.R.Acad.Sci.Paris},
volume={298},
date={1984},
pages={353-356},
}

\bib{Ca}{article}{
author={Carlsson G.},
title={Equivariant stable homotopy theory and related areas},
journal={Homology,Homotopy Appl.},
volume={3(2)},
date={2001},

}

\bib{CaCeGr}{article}{
author={Carrasco P.},
author={Cegarra A.M.},
author={Grandjean A.R.},
title={(Co)homology of crossed modules},
journal={J.Pure Appl.Algebra},
volume={168(2-3)},
date={2002},
pages={147-176},
}


\bib{CeGaOr}{article}{
author={Cegarra A.M.},
author={Garcia-Calcines J.M.}
author={Ortega J.A.}
title={On graded categorical groups and equivariant group extensions},
journal={Canad.J.Math.}
volume={54(3)},
date={2002},
pages={970-997},
}

\bib{CeGaOr}{article}{
author={Cegarra A.M.},
author={Garcia-Calcines J.M.}
author={Ortega J.A.}
title={Cohomology of groups with operators},
journal={Homology,Homotopy Appl.},
volume={4(1)},
date={2002},
pages={1-23},
}

\bib{CeGa}{article}{
author={Cegarra A.M.},
author={Garzon A.R.}
title={Equivariant group cohomology and Brauer group},
journal={Bull.Belg.Math.Soc.},
volume={10(3)},
date={2003},
pages={451-459},
}

\bib{CeIn}{article}{
author={Cegarra A.M.}
author={Inassaridze, H.},
title={Homology of groups with operators},
journal={Inter.Math.J.},
volume={5(1)},
date={2004},
pages={29-48},
}

\bib{DoLa}{article}{
author={Donadze G.},
author={Ladra M.},
title={More on the Hochschild and cyclic homology of crossed modules},
journal={Comm.Algebra},
volume={39(11)},
date={2011},
pages={4447-4460},
}

\bib{DoLi}{article}{
author={Donadze G.},
author={Van Der Linden T.},
title={A comonadic interpretation of Baues-Ellis homology of crossed modules},
journal={J.Homotopy,Rel,Str.},
volume={14},
date={2019},
pages={625-646},

}

\bib{EiGa}{article}{
author={Eilenberg S.},
author={Ganea T.},
title={On the Lusternik-Schnirelmann category of abstract groups},
journal={Ann. of Math.},
volume={65(2)},
date={1957},
pages={517-518}
}

\bib{El}{article}{
author={Ellis G.J.},
title={Homology of 2-types},
journal={J.London Math.Soc.},
volume={46(2)},
date={1992},
pages={1-27},
}
\bib{FiHaMa}{article}{
author={Fiedorowicz Z.},
author={Hauschild H.}
author={May P.}
title={Equivariant K-theory},
journal={Lecture Notes in Math. Springer}
date={1982},
pages={23-80},
}


\bib{Gi}{article}{
author={Gilbert N.D.},
title={The low-dimensional homology of crossed modules},
journal={Homology, Homotopy Appl.},
volume={2(1)},
date={2000},
pages={41-50},
}

\bib{GrLaPi}{article}{
author={Grandjean A.R.},
author={Ladra M.},
author={Pirashvili T.},
title={CCG-homology of crossed modules via classifying spaces},
journal={J.Algebra},
volume={229(2)},
date={2000},
pages={660-665},
}

\bib{GrMePa}{article}{
author={Grant M.},
author={Meir E.},
author={Patchkoria I.},
title={Equivariant dimensions of groups with operators},
journal= {Groups, Geometry,and Dynamics, accepted / in press, arXiv:1912.01692v3}
date={2020},

}

\bib{Gr}{book}{
author={Greither C.},
title={Cyclic Galois extensions of commutative rings},
Publisher={Lecture Notes in Mathematics,154,Springer},
date={1992},
}



\bib{Gu}{article}{
author={Guin D.},
title={Cohomologie et homologie non abéliennes des groupes},
journal={J.Pure Appl.Algebra},
volume={50},
date={1088},
pages={109-137},
}

\bib{Ho}{article}{
author={Hochschild G.},
title={Relative homological algebra},
journal={Trans.Amer.Math.Soc.},
volume={82},
date={1956},
pages={246-269},
}

\bib{Hu}{article}{
author={Huebschmann J.},
title={Crossed n-fold extensions of groups and cohomology},
journal={Commentarii Mathematici Helvetici},
volume={55},
date={1980},
pages={302-313},
}

\bib{In}{article}{
author={Inassaridze H.},
title={Homotopy of pseudosimplicial groups, non-abelian derived functors and algebraic K-theory},
journal={Math.USSR Sb },
volume={98(3)},
date={1975},
pages={339-362},
}

\bib{In}{book}{
author={Inassaridze H.},
title={Algebraic K-theory},
Publisher={Kluwer Academic Publishers, Dordrecht},
date={1995},
}

\bib{In}{book}{
author={Inassaridze H.},
title={Non-Abelian Homological Algebra and its Applications},
Publisher={Kluwer Academic, Dordrecht},
date={1997},

}

\bib{In}{article}{
author={Inassaridze H.},
title={Equivariant homology and cohomology of groups},
journal={Topology and its Applications },
volume={153},
date={2005},
pages={66-89},
}

\bib{InKh}{article}{
author={Inassaridze H.},
author={Khmaladze E.}
title={Hopf formulas for equivariant integral homology of groups},
journal={Proc.Amer.Math.Soc.},
volume={138(9)},
date={2010},
pages={3037-3046},
}

\bib{Ja}{article}{
author={Janelidze G.},
title={On abelian extensions of commutative rings},
journal={Bull.Georgian Acad.Sci.},
volume={108(2)},
date={1982},
pages={477-480}
}

\bib{Ke}{article}{
author={Keune F.},
title={The relativization of $K_{2}$},
journal={J.Algebra},
volume={54},
date={1978},
pages={159-177},
}

\bib{Ku}{article}{
author={Kuku A.},
title={Equivariant K-theory and the cohomology of profinite groups},
journal={Lecture Notes in Math. Springer},
volume={342},
date={1984},
pages={235-244},
}


\bib{LaGr}{article}{
author={Ladra M.},
author={Grandjean A.R.}
title={Crossed modules and homology},
journal={J.Pure Appl.Algebra},
volume={95},
date={1994},
pages={41-55},

}

\bib{Lo}{article}{
author={Loday J.-L.},
title={Cohomologie et groupe de Steinberg relatifs},
journal={J.Algebra},
volume={54},
date={1978},
pages={178-202},
}


\bib{Ma}{article}{
author={Matsumoto H.},
title={Sur les sous groupes arithmétiques des groupes semi-simples déployés},
journal={Ann.Sci.Norm.Sup.},
volume={2(4)},
date={1969},
pages={1-62},
}


\bib{Ma}{book}{
author={May J.P.},
title={Simplicial objects in algebraic topology},
Publisher={Van Nostrand Math.Studies, No11, Van Nostrand },
date={1967},
}

\bib{MiLi}{article}{
author={Di Micco D.},
author={Van Der Linden T.},
title={Universal central extensions of internal crossed modules via non-abelian tensor product},
journal={Appl.Cat.Structures},
volume={28}
date={2020},
}

\bib{Mi}{book}{
author={Milnor J.},
title={Introduction to Algebraic K-theory},
Publisher={Annals of Mathematic Studies 72, Princeton University Press, NJ},
date={1971},

}

\bib{MuWo}{article}{
author={Mueller L.},
author={Woike L.},
title={Equivariant higher Hochschild homology and topological field theories},
journal={Homology,Homotopy Appl.},
volume={22(1)},
date={2020},
pages={27-54}
}

\bib{PaYa}{article}{
author={Pamuk S.},
author={Yalcn E.},
title={Relative group cohomology and orbit category},
journal={Comm.Algebra},
volume={42(2)},
date={2014},
pages={3220-3243}
}

\bib{Pa}{article}{
author={Paoli S.},
title={(Co)homology of crossed modules with coefficients in a $\pi_{1}$-module},
journal={Homology, Homotopy Appl.},
volume={5(1)}
date={2003},
pages={261-296},
}

\bib{Ph}{article}{
author={Philips N.C.},
title={Equivariant K-theory and freeness of group actions on $C^{*}$-algebras},
journal={Lecture Notes in Math. Springer},
volume={1274},
date={1987},
}
\bib{Qu}{article}{
author={Quillen D.},
title={Higher algebraic K-theory},
journal={Lecture Notes in Math. Springer},
volume={341},
date={1973},
pages={77-139},
}

\bib{St}{article}{
author={Stallings J.R.},
title={On torsion-free groups with infinitely many ends},
journal={Ann. of Math.},
volume={88(2)},
date={1968},
pages={312-334},
}

\bib{St}{article}{
author={Stammbach U.},
title={Homology in Group Theory},
journal={Lecture Notes in Math.,Springer},
volume={359},
date={1973},

}

\bib{Sw}{article}{
author={Swan R.G.},
title={Groups of cohomological dimension one},
journal={J.Algebra},
volume={12},
date={1969},
pages={585-610}
}

\bib{Sw}{article}{
author={Swan R.G.},
title={Excision in algebraic K-theory},
journal={J.Pure Appl.Algebra},
volume={1}
date={1971},
pages={221-252}
}

\bib{TiVo}{article}{
author={Tierney M.},
author={Vogel W.}
title={Simplicial resolutions and derived functors},
journal={Math.Zeit.},
volume={111}
date={1969},
pages={1-44}
}

\bib{Wh}{article}{
author={Whitehead J.H.C.},
title={On group extensions with operators},
journal={Quart.J.Math.Oxford},
volume={1(2)}
date={1950},
pages={219-228}
}


\bib{Za}{article}{
author={Zavalo S.T.},
title={S-free operator groups I},
journal={Ukr. Mat. Zh.},
volume={16(5)}
date={1964},
pages={593-602}
}

\bib{Za}{article}{
author={Zavalo S.T.},
title={S-free operator groups II},
journal={Ukr. Mat. Zh.},
volume={16(6)}
date={1964},
pages={730-751}
}

\end{biblist}
\end{bibdiv}
\end{document}